\documentclass [12pt]{amsart}

\usepackage [T1] {fontenc}
\usepackage {amssymb}
\usepackage {amsmath}
\usepackage {amsthm}
\usepackage{tikz}
\usepackage{hyperref}
\usepackage{url}
\usepackage{xcolor}

\usepackage[titletoc,toc,title]{appendix}

\usetikzlibrary{arrows,decorations.pathmorphing,backgrounds,positioning,fit,petri}

\usetikzlibrary{matrix}
\usepackage{geometry}
\usepackage{enumitem}

\DeclareMathOperator {\td} {td}

\DeclareMathOperator {\Th} {Th}
\DeclareMathOperator {\ldim} {ldim}
\DeclareMathOperator {\tp} {tp}
\DeclareMathOperator {\etp} {etp}
\DeclareMathOperator {\qftp} {qftp}
\DeclareMathOperator {\ord} {ord}

\DeclareMathOperator {\cl} {cl}
\DeclareMathOperator {\im} {Im}
\DeclareMathOperator {\pr} {pr}

\DeclareMathOperator {\ACF} {ACF}

\DeclareMathOperator {\SL} {SL}

\DeclareMathOperator {\GL} {GL}
\DeclareMathOperator {\codim} {codim}

\DeclareMathOperator {\Loc} {Loc}
\DeclareMathOperator {\Exp} {Exp}

\DeclareMathOperator {\alg} {alg}

\DeclareMathOperator {\rank} {rank}

\newcommand{\D}{\mathrm{D}}

\theoremstyle {definition}
\newtheorem {definition}{Definition} [section]
\newtheorem {remark} [definition] {Remark}

\newtheorem{example} [definition] {Example}
\newtheorem* {claim} {Claim}
\newtheorem* {notation} {Notation}
\newtheorem {assumption} [definition] {Assumption}

\theoremstyle {plain}
\newtheorem {lemma} [definition] {Lemma}

\newtheorem {theorem} [definition] {Theorem}
\newtheorem {proposition} [definition] {Proposition}
\newtheorem {corollary} [definition] {Corollary}
\newtheorem {conjecture} [definition] {Conjecture}
\newtheorem {question} [definition] {Question}

\theoremstyle {remark}

\makeatletter
\newcommand {\forksym} {\raise0.2ex\hbox{\ooalign{\hidewidth$\vert$\hidewidth\cr\raise-0.9ex\hbox{$\smile$}}}}
\def\@forksym@#1#2{\mathrel{\mathop{\forksym}\displaylimits_{#2}}}
\def\forkind{\@ifnextchar_{\@forksym@}{\forksym}}
\newcommand {\nforkind} {\not \forkind}
\makeatother

\makeatletter
\def\@tocline#1#2#3#4#5#6#7{\relax
  \ifnum #1>\c@tocdepth 
  \else
    \par \addpenalty\@secpenalty\addvspace{#2}%
    \begingroup \hyphenpenalty\@M
    \@ifempty{#4}{%
      \@tempdima\csname r@tocindent\number#1\endcsname\relax
    }{%
      \@tempdima#4\relax
    }%
    \parindent\z@ \leftskip#3\relax \advance\leftskip\@tempdima\relax
    \rightskip\@pnumwidth plus4em \parfillskip-\@pnumwidth
    #5\leavevmode\hskip-\@tempdima
      \ifcase #1
       \or\or \hskip 1em \or \hskip 2em \else \hskip 3em \fi%
      #6\nobreak\relax
    \dotfill\hbox to\@pnumwidth{\@tocpagenum{#7}}\par
    \nobreak
    \endgroup
  \fi}
\makeatother

\calclayout

\begin {document}

\title[Predimensions in Differential Fields]{Adequate Predimension Inequalities in Differential Fields}

\author{Vahagn Aslanyan}

\thanks{Most of this work was done while the author was a DPhil student at the University of Oxford. Some revisions where done at Carnegie Mellon University and the University of East Anglia.\\ Supported by the University of Oxford Dulverton Scholarship, Luys Scholarship and AGBU UK Scholarship (at the University of Oxford) and EPSRC grant EP/S017313/1 (at the University of East Anglia).}
\address{{Department of Mathematical Sciences, CMU, Pittsburgh, PA 15213, USA}}

\address{{Institute of Mathematics of the National Academy of Sciences, Yerevan 0019, Armenia}}

\address{\texttt{Current address:} School of Mathematics, University of East Anglia, Norwich, NR4 7TJ, UK}
\email{V.Aslanyan@uea.ac.uk}




\vspace*{-1.5cm}

\maketitle

\begin{abstract}
In this paper we study predimension inequalities in differential fields and define what it means for such an inequality to be \emph{adequate}. Adequacy was informally introduced by Zilber, and here we give a precise definition in a quite general context. We also discuss the connection of this problem to definability of derivations in the reducts of differentially closed fields. The Ax-Schanuel inequality for the exponential differential equation (proved by Ax) and its analogue for the differential equation of the $j$-function (established by Pila and Tsimerman) are our main examples of predimensions. We carry out a Hrushovski construction with the latter predimension and obtain a natural candidate for the first-order theory of the differential equation of the $j$-function. It is analogous to Kirby's axiomatisation of the theory of the exponential differential equation (which in turn is based on the axioms of Zilber's pseudo-exponentiation), although there are many significant differences. In joint work with Sebastian Eterovi\'c and Jonathan Kirby we have recently proven that the axiomatisation obtained in this paper is indeed an axiomatisation of the theory of the differential equation of the $j$-function, that is, the Ax-Schanuel inequality for the $j$-function is adequate.
\end{abstract}

\tableofcontents


\section{Introduction}

In this paper we study predimension inequalities in differential fields and, in particular, \emph{adequacy} of such inequalities. The notion of a predimension was introduced by Hrushovski in \cite{Hru} where he used an amalgamation-with-predimension technique (which is a variation of Fra\"iss\'e's amalgamation construction) to construct counterexamples to Zilber's Trichotomy Conjecture. Given a predimension on a class of structures satisfying certain properties, one can carry out an amalgamation-with-predimension construction \`a la Hrushovski and obtain a limit structure which is homogeneous and existentially closed with respect to embeddings respecting the predimension (known as \emph{strong} embeddings). Depending on the predimension, this limit structure can have various model theoretic properties and may or may not be a ``natural'' mathematical structure. In a sense, adequacy means that the construction yields a natural mathematical structure whose properties are usually the main object of one's interest. We illustrate this idea on some well-known examples below postponing the definition of adequacy to \S \ref{predimensions-and-Hrush}.

In the early 2000's Zilber explored complex exponentiation from a model theoretic perspective and observed that Schanuel's conjecture plays a key role in understanding its model theoretic properties. The conjecture states that for any $\mathbb{Q}$-linearly independent complex numbers $z_1,\ldots,z_n$ the transcendence degree of $z_1,\ldots,z_n,e^{z_1},\ldots,e^{z_n}$ over $\mathbb{Q}$ is at least $n$ (see \cite[p. 30]{Lang-tr}). This can be written in the following equivalent form: for any $z_1,\ldots,z_n \in \mathbb{C}$ the inequality
\begin{equation*}\label{Schanuel-predimension}
\delta(\bar{z}) := \td_{\mathbb{Q}}\mathbb{Q}(\bar{z},\exp(\bar{z}))-\ldim_{\mathbb{Q}}(\bar{z}) \geq 0
\end{equation*}
holds, where $\td$ and $\ldim$ stand for transcendence degree and linear dimension respectively. Zilber observed that $\delta$ is a predimension function, and Schanuel's conjecture states that it is non-negative.\footnote{Although we will present a general axiomatic definition of predimensions in Section \ref{predimensions-and-Hrush}, our examples will look like the Schanuel predimension. So the reader can think of predimensions as combinations of dimension functions such as the transcendence degree of a field, the linear dimension of a vector space and the cardinality of a set.} Using Hrushovski's aforementioned technique, Zilber constructed algebraically closed fields of characteristic $0$ equipped with a unary function, called \emph{pseudo-exponentiation}, which mimics the basic properties of complex exponentiation and, more importantly, satisfies the analogue of Schanuel's conjecture. It is also \emph{Exponentially Algebraically Closed} which means roughly that a system of equations has a solution unless having a solution contradicts Schanuel's conjecture. Zilber conjectured that the unique pseudo-exponential field of cardinality $2^{\aleph_0}$ is isomorphic to $\mathbb{C}_{\exp} := (\mathbb{C}; +, \cdot, \exp )$. In this setting the ``natural mathematical structure'' is $\mathbb{C}_{\exp}$, and Zilber's conjecture states the adequacy of the Schanuel predimension (assuming Schanuel's conjecture).
See \cite{Zilb-pseudoexp} for details (pseudo-exponentiation will also be briefly discussed in \S \ref{Predim-examples}).

Although Schanuel's conjecture is considered out of reach, Ax proved its differential/functional analogue in \cite{Ax} which is now known as the Ax-Schanuel theorem. It is a transcendence inequality for solutions of the exponential differential equation in a differential field, and gives a predimension inequality on the class of exponential reducts of differential fields, that is, reducts of the form $\mathcal{K}_{\Exp} := (K; +, \cdot, \Exp)$ where $\mathcal{K} := (K; +, \cdot, ')$ is a differential field and $\Exp(x,y)$ is given by the formula $y' = y x' $. More precisely, given $(x_i, y_i) \in \Exp(K),~ i=1,\ldots,n,$ the predimension is given by 
\[
\delta(\bar{x}):= \td_CC(\bar{x},\bar{y}) - \ldim_{\mathbb{Q}}(x_1',\ldots,x_n'),
\]
where $C$ is the field of constants of $\mathcal{K}$. It is easy to see that $\delta$ does not depend on the choice of $y_i$'s. Then the Ax-Schanuel theorem states that $\forall \bar{x} (\delta(\bar{x})\geq 0)$ and equality holds if and only if $\bar{x}\in C^n$.

It is often possible to give an axiomatisation of the first-order theory of the Hrushovski-Fra\"iss\'e limit $\mathcal{U}$ and establish certain model theoretic properties such as stability. For the exponential Ax-Schanuel inequality this was done by Kirby in \cite{Kirby-semiab}, following Zilber's work on pseudo-exponentiation. In that case, the limit $\mathcal{U}$ is a structure in the language $\mathfrak{L}_{\Exp} := \{ +, \cdot, \Exp \}$. Kirby obtained an axiomatisation of $\Th(\mathcal{U})$ which consists of several basic axioms and, most importantly, the Ax-Schanuel inequality and a matching \emph{Existential Closedness} property, or EC for short. The latter states that certain systems of equations in the language $\mathfrak{L}_{\Exp}$ have solutions in $\mathcal{U}$, and is equivalent to $\mathcal{U}$ being existentially closed in strong extensions. Note that EC is the differential analogue of Zilber's \emph{Exponential Algebraic Closedness} axiom for pseudo-exponentiation. Further, Kirby showed that $\mathcal{U}$ is isomorphic to the reduct $\mathcal{F}_{\Exp}$ of the countable saturated differentially closed field $\mathcal{F}$. This means that the Ax-Schanuel inequality for $\Exp$ is adequate. In particular, $\mathcal{F}_{\Exp}$ satisfies EC, that is, certain systems of exponential differential equations have solutions in differentially closed fields (see also \cite[\S 4]{Aslanyan-Eterovic-Kirby-Diff-EC-j}). Actually this fact is equivalent to adequacy, so the latter can be understood as a property on solvability of some systems of equations. We will elaborate on this in \S \ref{Predim-examples}.

Adequacy was informally introduced by Zilber during numerous discussions around these topics and, in particular, the term ``adequate'' is due to him. It was also mentioned in the papers \cite{Aslanyan-def-deriv} and \cite{Aslanyan-linear}, and the latter even contains a definition of adequacy in the differential context. In this paper, we define adequacy more formally and more generally and pose the following question.

\begin{question}\label{main-question-0}
Which differential equations satisfy an adequate predimension inequality?
\end{question}

We will often refer to predimension inequalities in differential fields as Ax-Schanuel type inequalities. It is well known that this kind of inequalities (even without adequacy) have important applications in transcendence theory (see \cite{Bays-Kirby-Wilkie,Kirby-Schanuel,Eterovic-Schan-for-j}) and diophantine geometry (see \cite{Pila-functional-transcendence,Zilb-exp-sum-published,Kirby-semiab,Pila-Tsim-Ax-j,Habegger-Pila-o-min-certain,Zannier-book-unlikely,Daw-Ren,Aslanyan-weakMZPD,Aslanyan-remarks-atyp}). Furthermore, adequacy gives a good understanding of the model theoretic and geometric properties of the differential equation under consideration. In particular, as in the case of the exponential differential equation, adequacy is normally equivalent to an Existential Closedness statement which states roughly that systems of equations in the appropriate language always have solutions unless having a solution contradicts the Ax-Schanuel inequality. Two more applications of adequacy will be mentioned below.

In \cite{Aslanyan-def-deriv} we explored the question of definability of derivations in reducts of differentially closed fields containing the field structure. We also explained there how that question is related to adequacy of Ax-Schanuel type inequalities. In this paper we will use some results from \cite{Aslanyan-def-deriv} to prove formally that under certain assumptions, if a derivation is definable in a reduct of a differentially closed field, then there can be no non-trivial adequate predimension inequality in that reduct. In \cite{Aslanyan-linear} we generalised Ax-Schanuel to linear differential equations with constant coefficients of higher order, and proved their adequacy.

Another example of a predimension inequality is the modular analogue of the Ax-Schanuel theorem, established by Pila and Tsimerman in \cite{Pila-Tsim-Ax-j}, that will be discussed in detail in \S \ref{j-chapter}. One of our main goals in this paper is to carry out a Hrushovski construction with that predimension and establish some model theoretic properties of the limit structure $\mathcal{U}$. In particular, we obtain an axiomatisation of $\Th(\mathcal{U})$, formulate an EC property and prove, among other things, the following theorem (see Theorems \ref{thm: j axiomatisation}, \ref{thm: j main}, \ref{thm: j' axiomatisation}, \ref{thm: j' main} for more general statements).

\begin{theorem}\label{thm-intro}
The Ax-Schanuel inequality for $j$ is strongly adequate if and only if the appropriate reducts of differentially closed fields satisfy the EC property. 
\end{theorem}

The EC property is a first-order axiom scheme stating that \emph{broad} algebraic varieties contain points which are solutions of the differential equation of the $j$-function. Broadness is an algebraic property dictated by Ax-Schanuel stating that the variety is not too small. Thus, as it was explained above, EC can be understood as a statement about systems of equations having a solution unless the existence of a solution contradicts Ax-Schanuel. Broadness and other related notions will be defined later in the paper.

Note that even though Theorem \ref{thm-intro} is analogous to Kirby's results for the exponential differential equation, its proof significantly differs from Kirby's proof, and we also take this opportunity to fix some inaccuracies from \cite{Kirby-semiab}. Furthermore, the $j$-function satisfies an order $3$ differential equation so one can consider two versions of the problem, one with the $j$-function only and one with $j$ and its first two derivatives. The latter is considerably more involved, but we consider both versions.

It was recently proven in our joint work with Eterovi\'c and Kirby \cite{Aslanyan-Eterovic-Kirby-Diff-EC-j} that EC holds in differentially closed fields. Combining that with Theorem \ref{thm-intro} we get adequacy.

\begin{theorem}
The Ax-Schanuel inequality for the $j$-function is adequate.\footnote{The current paper was initially written before \cite{Aslanyan-Eterovic-Kirby-Diff-EC-j} and even before adequacy was proven. So in earlier versions of this paper adequacy was stated as a conjecture.}
\end{theorem}

Adequacy of the modular Ax-Schanuel inequality has several important applications including a Zilber-Pink type statement for the $j$-function and its derivatives \cite{Aslanyan-weakMZPD} and a characterisation of strongly minimal sets in $j$-reducts of differentially closed fields \cite{Aslanyan-SM-reducts}.

\addtocontents{toc}{\protect\setcounter{tocdepth}{1}}
\subsection*{Outline of the paper}

In Section \ref{predimensions-and-Hrush} we give a brief account of predimensions and Hrushovski style amalgamation-with-predimension constructions. In particular, we give a rigorous definition of adequacy of a predimension inequality. We also consider some examples and show how they fit with the presented approach. 

In Section \ref{connection-def-der} we discuss the connection of Question \ref{main-question-0} to the question of definability of a derivation in reducts of differentially closed fields. We show essentially that if a derivation is definable in a reduct of a differentially closed field then the reduct cannot satisfy any non-trivial adequate predimension inequality.

Sections \ref{j-chapter} and \ref{j-general-case} are devoted to the differential equation of the $j$-function. Starting with the Ax-Schanuel inequality for $j$, we show that the class of models of a certain theory (which is essentially the universal theory of reducts of differential fields with a relation for the equation of $j$)  has the strong amalgamation property. Then we construct the Hrushovski-Fra\"iss\'e limit and give an axiomatisation of its first-order theory. The given axiomatisation is a natural candidate for the theory of the differential equation of the $j$-function, and in recent joint work with Sebastian Eterovi\'c and Jonathan Kirby we proved that it is indeed an axiomatisation of that theory. In other words, we proved that the Ax-Schanuel inequality for the $j$-function is adequate (see \cite{Aslanyan-Eterovic-Kirby-Diff-EC-j}).

Most of this work forms part of the author's DPhil thesis \cite{Aslanyan-thesis}.

\subsection*{Notation and conventions}
\setcounter{equation}{0}

Here we fix some notation that will be used throughout the paper. 
\begin{itemize}[leftmargin = 0.5cm]


\item The length of a tuple $\bar{a}$ will be denoted by $| \bar{a}|$. For a set $A$ and a tuple $\bar{a}$ we will sometimes write $\bar{a} \in A$ or $\bar{a} \subseteq A$ and mean that all coordinates of $\bar{a}$ are in $A$, i.e. $\bar{a} \in A^{| \bar{a} |}$.

\item For two sets $X,Y$ the notation $X \subseteq_{fin} Y$ means $X$ is a finite subset of $Y$. The union $X \cup Y$ will sometimes be written as $XY$. The power set of $X$ is denoted by $\mathfrak{P}(X)$.

\item 
For the linear dimension of a vector space $V$ over a field $K$ we use the shorthand $\ldim_KV$.

\item All fields considered in this work will be of characteristic zero. The algebraic closure of a field is denoted by $K^{\alg}$.

\item By an irreducible variety we always mean absolutely irreducible.

\item If $K \subseteq F$ are fields, the transcendence degree of $F$ over $K$ will be denoted by $\td_KF$ or $\td(F/K)$. When we work in an ambient algebraically closed field $F$ and $V$ is a variety defined over $F$, we will normally identify $V$ with the set of its $F$-points $V(F)$. The algebraic locus (Zariski closure) of a tuple $\bar{a} \in F$ over $K$ will be denoted by $\Loc_K (\bar{a})$ or $\Loc(\bar{a}/K)$ (and identified with the set of its $F$-points).

\item If $F$ is a differential field then for a non-constant element $x \in F$ the \emph{differentiation with respect to} $x$ is a derivation $\partial_x: F \rightarrow F$ defined by $y\mapsto \frac{y'}{x'}$, where $'$ is the derivation of $F$.

\item For a differential field $F$ the differential ring of differential polynomials over $F$ of variables $\bar{X}$ is denoted by $F\{ \bar{X} \}$.



\item If $\mathcal{M}$ is a structure and $\bar{a} \in M^n$ is a finite tuple, then the complete type of $\bar{a}$ in $\mathcal{M}$ over a parameter set $A \subseteq M$ will be denoted by $\tp^{\mathcal{M}}(\bar{a}/A)$ while $\qftp^{\mathcal{M}}(\bar{a}/A)$ stands for the quantifier-free type. We often omit the superscript $\mathcal{M}$ if the ambient model is clear.


\end{itemize}
\addtocontents{toc}{\protect\setcounter{tocdepth}{2}}

\section{Predimensions and Hrushovski constructions}\label{predimensions-and-Hrush}

In this section we present the appropriate definitions of predimensions and strong embeddings and observe several standard facts about them. Then we give a brief account of Hrushovski's amalgamation-with-predimension construction. It is the uncollapsed version of a full Hrushovski construction \cite{Hru,Hru-acf}. This will be used to define \emph{adequacy} of a predimension inequality. The Ax-Schanuel inequality for the exponential differential equation (Ax \cite{Ax}) and its analogue for the differential equation of the modular $j$-function (Pila-Tsimerman \cite{Pila-Tsim-Ax-j}) are our main examples. We will observe the close relationship between triviality of an adequate predimension inequality and model completeness of the corresponding strong Fra\"iss\'e limit. 


We mainly follow Wagner \cite{Wagner} and Baldwin \cite{Baldwin-rank} in defining predimensions and related notions. They give an axiomatic approach to Hrushovski constructions. Wagner works in a relational language, while Baldwin's setting does not have this restriction. We need that generality since we always have a field structure in our examples. Note that Baldwin imposes stronger definability conditions for the predimension than we do. The reason is that the Ax-Schanuel predimension does not satisfy his definability axioms. Our approach is motivated by Kirby's analysis of the exponential differential equations \cite{Kirby-semiab} and Zilber's approach to complex exponentiation and Schanuel's conjecture \cite{Zilb-pseudoexp}.

Hrushovski invented the aforementioned constructions in order to produce structures with ``exotic'' geometry and refute some conjectures on categorical theories and answer some questions. Most notably, he refuted Zilber's Trichotomy Conjecture \cite{Zilber-sm,Zilber-ICM} stating that any uncountably categorical and non-locally modular theory is bi-interpretable with an algebraically closed field, and Lachlan's conjecture \cite{Lachlan-conj} stating that any stable $\aleph_0$-categorical theory is totally transcendental. Later on Hrushovski's techniques were adapted and used in various settings to construct interesting structures. The reader is referred to \cite{Hru-acf,Hru,Wagner,Wagner-Hrush,Baldwin-rank,Baldwin-Holland,Zilb-pseudoexp} for details on Hrushovski constructions and examples of ``exotic'' structures (theories) that can be obtained by such constructions.

Most of the concepts and results of this section are quite standard and well-known, and we present them here and often prove standard facts to fix the notation and make the paper self-contained.

\subsection{Predimensions}\label{predimensions}



Let $\mathfrak{L}$ be a countable language and $\mathfrak{C}$ be a collection of $\mathfrak{L}$-structures closed under intersections. This can be understood in a category theoretic sense, but for us it will be enough to assume that if $A_i \in \mathfrak{C},~ i \in I$, are substructures of some $A \in \mathfrak{C}$ then $\bigcap_{i\in I} A_i \in \mathfrak{C}$. We will also assume that $\mathfrak{C}$ has the joint embedding property, i.e. for any $A, B \in \mathfrak{C}$ there is $C \in \mathfrak{C}$ such that $A$ and $B$ can be embedded into $C$. Assume further that $\mathfrak{C}$ contains a smallest structure $S \in \mathfrak{C}$, that is, for every $A \in \mathfrak{C}$ we have $S \subseteq A$. Finally, $\mathfrak{C}$ is assumed to be closed under isomorphism over $S$, that is, if $A \subseteq \mathfrak{C}$ and $B \supseteq S$ is a structure and $f: A \rightarrow B$ is an isomorphism which fixes $S$ pointwise, then $B \in \mathfrak{C}$.

\begin{definition}
For $B \in \mathfrak{C}$ and $X \subseteq B$ the $\mathfrak{C}$-\emph{closure} of $X$ inside $B$ (or the $\mathfrak{C}$-substructure of $B$ generated by $X$) is the structure\footnote{In the differential setting the notation $\langle A \rangle$ is used to denote the differential subfield generated by a set $A$. The meaning of this notation will be clear from the context. In particular it is used only for $\mathfrak{C}$-closure in this section.}
$$ \langle X \rangle_B := \bigcap_{A \in \mathfrak{C}: X \subseteq A \subseteq B}A.$$
A structure $A \in \mathfrak{C}$ is \emph{finitely generated} if $A=\langle X \rangle_A$ for some finite $X \subseteq A$.
{The collection of all finitely generated structures from $\mathfrak{C}$ will be denoted by $\mathfrak{C}_{f.g.}$.} For $A, B \in \mathfrak{C}$ by $A \subseteq_{f.g.} B$ we mean $A$ is a finitely generated substructure of $B$.
\end{definition}

Note that in general finitely generated in this sense is different from being finitely generated as a structure. We will assume however that finitely generated structures are countable. A substructure of a finitely generated structure may not be finitely generated. However, in our examples this does not happen. So we make the following assumption.

\begin{assumption}\label{fg}
Assume $\mathfrak{C}$ satisfies the following condition.
\begin{itemize}[leftmargin = 1cm]
\item[FG] If $A \in \mathfrak{C}, B \in \mathfrak{C}_{f.g.}$ with $A \subseteq B$ then $A \in \mathfrak{C}_{f.g.}$.
\end{itemize}
\end{assumption} 

{Since $S$ is the smallest structure in $\mathfrak{C}$, it is in fact generated by the empty set, i.e. $S = \langle \emptyset \rangle$. So, by abuse of notation, we will normally write $\emptyset$ instead of $S$.}

When we have two structures $A,B \in \mathfrak{C}$ we would like to have a notion of a structure generated by $A$ and $B$. However, this cannot be well-defined without embedding $A$ and $B$ into a bigger $C$. Given such a common extension $C$, we will denote $AB_C:= \langle A\cup B \rangle_C$. Often we will drop the subscript $C$ meaning that our statement holds for every common extension $C$ (or it is obvious in which common extension we work). This remark is valid also when we write $A \cap B$ which should be understood as the intersection of $A$ and $B$ after identifying them with their images in a common extension.

\begin{definition}
A \emph{predimension} on $\mathfrak{C}_{f.g.}$ is a function $\delta: \mathfrak{C}_{f.g.} \rightarrow \mathbb{Z}$ with the following properties:
\begin{enumerate}[leftmargin = 1cm]
\item[P1] $\delta(\emptyset) = 0$,
\item[P2] If $A, B \in \mathfrak{C}_{f.g.}$ with $A \cong B$ then  $\delta(A)= \delta(B)$,
\item[P3] (Submodularity) For all $A, B \in \mathfrak{C}_{f.g.}$ and $C \in \mathfrak{C}$ with $A, B \subseteq C$ we have  
\begin{equation}\label{submodularity}
\delta(AB)+ \delta(A\cap B) \leq \delta(A) + \delta (B).
\end{equation}
\end{enumerate}
\end{definition}

If, in addition, such a function is monotonic, i.e. $A \subseteq B \Rightarrow \delta (A) \leq \delta(B)$, and hence takes only non-negative values, then  $\delta$ is called a \emph{dimension}. 

\begin{remark}
If equality holds in \eqref{submodularity} for all $A, B\in \mathfrak{C}_{f.g.}$ then $\delta$ is said to be \emph{modular}.
\end{remark}

\begin{definition}
Given a predimension $\delta$, for a finite subset $X \subseteq_{fin} A \in \mathfrak{C}$ one defines
$$\delta_A(X):=\delta(\langle X \rangle_A).$$
\end{definition}

The following is Hrushovski's \emph{ab initio} example from \cite{Hru}.

\begin{example}
Let $\mathfrak{C}$ be the class of all structures $(M;R)$ in a language $\mathfrak{L}=\{ R\}$ consisting of one ternary relation $R$. Then $\mathfrak{C}_{f.g.}$ is the collection of all finite $\mathfrak{L}$-structures. For $A \in \mathfrak{C}_{f.g.}$ define
$$\delta(A):=|A|-|R(A)|.$$
Then $\delta$ is a predimension.
\end{example}

Other examples of predimensions, which are more relevant to our work, will be given in Section \ref{Predim-examples}.


{Now we define the relative predimension of two structures, which depends on a common extension of those structures (so we work in such a common extension without explicitly mentioning it). }

\begin{definition}
Let $\delta$ be a predimension on $\mathfrak{C}_{f.g.}$. The \emph{relative predimension} is defined as follows.
\begin{itemize}[leftmargin = 0.5cm]
\item For $A, B \in \mathfrak{C}_{f.g.}$ define $\delta(A/B):=\delta(AB) - \delta(B)$. 

\item For $X \in \mathfrak{C}_{f.g.}$ and $A \in \mathfrak{C}$ define $\delta(X/A) \geq k$ for an integer $k$ if for all $Y\subseteq_{f.g.}A$ there is $Y \subseteq Y' \subseteq_{f.g.}A$ such that $\delta(X/Y') \geq k$. Define $\delta(X/A) = k$ if $\delta(X/A) \geq k$ and $\delta(X/A) \ngeq k+1$.
\end{itemize}
\end{definition}

In the next definition $B$ is the ambient structure that we work in.

\begin{definition}
Let $A \subseteq B \in \mathfrak{C}$. We say $A$ is \emph{strong} (or \emph{self-sufficient}) in $B$, denoted $A \leq B$, if for all $X \subseteq_{f.g.} B$ we have $\delta(X/A) \geq 0$. One also says $B$ is a \emph{strong extension} of $A$. An embedding $A \hookrightarrow B$ is strong if the image of $A$ is strong in $B$.
\end{definition}

It is easy to notice that the above definition will not change if we take a finite set $X$ instead of a finitely generated structure $X$.

\begin{lemma}
Let $A, B \in \mathfrak{C}$. Then $A \leq B$ if and only if for all $X \subseteq_{f.g.} B$ we have $\delta(X\cap A) \leq \delta(X)$.
\end{lemma}
\begin{proof}
Let $A \leq B$ and $X \subseteq_{f.g.} B$. Choose $Y=X \cap A \subseteq_{f.g.} A$. Then by definition there is $Y \subseteq Y' \subseteq_{f.g.} A$ such that $\delta(XY')\geq \delta(Y')$. Now by submodularity
$$\delta(X \cap A)= \delta(Y)=\delta(X\cap Y') \leq \delta(X)+\delta(Y')-\delta(XY')\leq \delta(X).$$

Conversely, assume the condition given in the lemma holds. We need to prove that $A \leq B$. Let $X \subseteq_{f.g.} B$ and $Y \subseteq_{f.g.} A$. Choose $Y'=XY \cap A \supseteq Y$. Then
$$\delta(XY')=\delta(XY) \geq \delta(XY \cap A) = \delta(Y')$$
as $XY \subseteq_{f.g.} B$.
\end{proof}

Observe that $\leq$ is transitive, i.e. if $A\leq B \leq C$ then $A\leq C$.

Now we assume $\mathfrak{C}_{f.g.}$ satisfies the Descending Chain Condition.

\begin{assumption}\label{DCC}
\begin{itemize}[leftmargin = 1.5cm]
    \item[]
    \item[DCC] $\mathfrak{C}_{f.g.}$ does not contain any infinite  strictly descending strong chain.
\end{itemize}

\end{assumption}

\begin{definition}
For $B \in \mathfrak{C}$ and $X \subseteq B$ we define the \emph{self-sufficient closure} of $X$ in $B$ by 
$$\lceil X \rceil_B := \bigcap_{A \in \mathfrak{C}: X \subseteq A \leq B}A.$$
\end{definition}

\begin{lemma}
For $B\in \mathfrak{C}$ and $X\subseteq B$ we have $\lceil X \rceil_B \leq B$.
\end{lemma}
\begin{proof}
It is easy to see that the intersection of finitely many strong substructures is strong. We need to prove that an arbitrary intersection of strong substructures is strong. Assume $A_i\leq B,~ i \in I$. Set $A:=\bigcap_{i\in I} A_i$. Pick $X\subseteq_{f.g.} B$ and denote $X_i:=X\cap A_i$. Then $X_i \leq X$ and
$$X\cap A = \bigcap_{i\in I} (X \cap A_i) = \bigcap_{i\in I} X_i.$$
By DCC there is a finite subset $I_0\subseteq I$ with $$\bigcap_{i\in I} X_i = \bigcap_{i\in I_0} X_i \leq X.$$ 
\end{proof}


\begin{lemma}\label{predim-type}
Let $M \in \mathfrak{C}$ be saturated. If $X,Y \subseteq_{fin} M$ (with some enumeration) have the same type in $M$  then $\langle X \rangle_M \cong \langle Y \rangle_M$ and $\lceil X \rceil_M \cong \lceil Y \rceil_M$ and hence $\delta_M(X)=\delta_M(Y)$.
\end{lemma}
\begin{proof}
Since $M$ is saturated and $\tp(X)=\tp(Y)$, there is an automorphism that sends $X$ to $Y$. Now the lemma follows from P2.
\end{proof}

From now on we assume $\delta(A) \geq 0$ for all $A \in \mathfrak{C}_{f.g.}$. In other words $\emptyset$ is strong in all structures of $\mathfrak{C}$. Instead of assuming this we could work with the subclass $\mathfrak{C}^0$ of all structures with non-negative predimension. However, we find it more convenient to assume $\delta$ is non-negative on $\mathfrak{C}_{f.g.}$ since it will be the case in our examples anyway.

\begin{lemma}\label{sscl-min}
If $B \in \mathfrak{C}$ and $X \subseteq_{f.g.} B$ then
\begin{itemize}[leftmargin = 0.5cm]
\item  $\lceil X\rceil_B$ is finitely generated, and

\item $\delta(\lceil X \rceil_B)=\min\{\delta(Y): X\subseteq Y \subseteq_{f.g.} B\}.$
\end{itemize}

\end{lemma}
\begin{proof}
Let $A\subseteq_{f.g.} B$ be such that $\delta(A)=\min \{ \delta(A'): X\subseteq A' \subseteq_{f.g.} B \}$. We claim that $A \leq B$. Indeed, for any $Y \subseteq_{f.g.} B$ we have
$$\delta(A\cap Y) \leq \delta(A)-\delta(AY)+\delta(Y)\leq \delta(Y).$$
Thus $A \leq B$ and hence $\lceil X \rceil_B$ is contained in finitely generated $A$ and so is finitely generated itself.

Further, $\lceil X \rceil_B \leq A$ so $\delta(\lceil X \rceil_B) \leq \delta(A)$. Now by minimality of $\delta(A)$ we conclude that $\delta(\lceil X \rceil_B) = \delta(A)$.
\end{proof}



A predimension gives rise to a dimension in the following way. 
\begin{definition}
For $X \subseteq_{f.g.} B$ define
$$d_B(X):= \min \{ \delta(Y): X \subseteq Y \subseteq_{f.g.} B \}=\delta(\lceil X\rceil_B).$$
For $X \subseteq_{fin} B$ set $d_B(X):=d_B(\langle X \rangle_B)$. 
\end{definition}

It is easy to verify that $d$ is a dimension function and therefore we have a natural pregeometry associated with $\delta$. More precisely, we define $\cl_B : \mathfrak{P}(B) \rightarrow \mathfrak{P}(B)$ by 
$$\cl_B(X)=\{ b \in B: d_B(b/X)=0 \}.$$
Then $(B,\cl_B)$ is a pregeometry and $d_B$ is its dimension function.

Self-sufficient embeddings can be defined in terms of $d$. Indeed, if $A \subseteq B$ then $A \leq B$ if and only if for any $X \subseteq_{fin} A$ one has $d_A(X)=d_B(X)$.

\begin{definition}\label{defin-predim-trivial-1}
A predimension $\delta$ is \emph{trivial} if all embeddings are strong. Equivalently, $\delta$ is trivial if it is monotonic and hence equal to the dimension associated with it.
\end{definition}

\begin{proposition}\label{elementary-strong}
Let $A, B \in \mathfrak{C}$ be saturated and $A \preceq B$. Then $A\leq B$.
\end{proposition}
\begin{proof}
If $A \nleq B$ then for some $X \subseteq_{f.g.} B$ one has $\delta(X/A)<0$. This means that there is $Y \subseteq_{f.g.} A$ such that for all $Y \subseteq Y' \subseteq_{f.g.} A$ we have $\delta(X/Y')<0$. Choose $Y'=\lceil Y \rceil_A$. The latter is finitely generated by Lemma \ref{sscl-min}. Suppose $X=\langle \bar{x} \rangle_B$ and $Y'=\langle \bar{y} \rangle_A$ for some finite tuples $\bar{x}$ and $\bar{y}$. Let $\bar{z}$ be a realisation of the type $\tp^B(\bar{x}/\bar{y})$ in $A$. If $Z=\langle \bar{z} \rangle_A$ then $\delta(Z/Y')<0$ (by Lemma \ref{predim-type}) which means $\delta(Y'Z)< \delta(Y')$ contradicting Lemma \ref{sscl-min}.
\end{proof}

\subsection{Amalgamation with predimension}\label{amalgamation}
Now we formulate conditions under which one can carry out an amalgamation-with-predimension construction. Let $\mathfrak{C}$ be as above and let $\delta$ be a non-negative predimension on $\mathfrak{C}_{f.g.}$.

\begin{definition}
The class $\mathfrak{C}$ is called a \emph{strong amalgamation class} if the following conditions hold.

\begin{itemize}[leftmargin = 1.5cm]
\item[C1] Every $A \in \mathfrak{C}_{f.g.}$ has at most countably many finitely generated strong extensions up to isomorphism.

\item[C2] $\mathfrak{C}$ is closed under unions of countable strong chains $A_0 \leq A_1 \leq \ldots$.

\item[SAP] $\mathfrak{C}_{f.g.}$ has the \emph{strong amalgamation property}, that is, for all $A_0, A_1, A_2 \in \mathfrak{C}_{f.g.}$ with $A_0 \leq A_i,~ i=1,2$, there is $B \in \mathfrak{C}_{f.g.}$ such that $A_1$ and $A_2$ are strongly embedded into $B$ and the corresponding diagram commutes.
\end{itemize}

\end{definition}

\begin{remark}
Since $\delta(A)\geq 0$ for all $A \in \mathfrak{C}_{f.g.}$, it follows that $\emptyset$ is strong in all finitely generated structures and hence the strong amalgamation property implies the \emph{strong joint embedding} property. 
\end{remark}

The following is a standard theorem that follows in particular from the category theoretic version of Fra\"iss\'e's amalgamation construction due to Droste and G\"obel \cite{Droste-Gobel} (see \cite{Kirby-semiab} for a nice exposition, without a proof though).

\begin{theorem}[Amalgamation theorem]
If $\mathfrak{C}$ is a strong amalgamation class then there is a unique (up to isomorphism) countable structure $U \in \mathfrak{C}$ with the following properties.
\begin{enumerate}[leftmargin = 1cm]
\item[\emph{U1}] $U$ is \emph{universal} with respect to strong embeddings, i.e. every countable $A \in \mathfrak{C}$ can be strongly embedded into $U$.
\item[\emph{U2}] $U$ is \emph{saturated} with respect to strong embeddings, i.e. for every $A, B \in \mathfrak{C}_{f.g.}$ with strong embeddings $A \hookrightarrow U$ and $A \hookrightarrow B$ there is  a strong embedding of $B$ into $U$ over $A$.
\end{enumerate}

Furthermore, any isomorphism between finitely generated strong substructures of $U$ can be extended to an automorphism of $U$.
\end{theorem}

This $U$ is called the \emph{generic model}, \emph{strong amalgam}, \emph{strong Fra\"iss\'e limit} or \emph{Fra\"iss\'e-Hrushovski limit} of $\mathfrak{C}_{f.g.}$. It has a natural pregeometry associated with the predimension function as described in the previous section. Note that U2 is normally known as the \emph{richness} property in literature (we used the terminology of \cite{Droste-Gobel} above).

\begin{remark}
Since we have assumed $\emptyset$ is strong in all structures from $\mathfrak{C}$, the property U2 implies U1. Indeed, for $A \in \mathfrak{C}_{f.g.}$ we have $\emptyset \leq A$ and $\emptyset \leq U$. Hence by U2 there is a strong embedding $A \hookrightarrow U$. Now since every countable structure in $\mathfrak{C}$ is the union of a strong chain of finitely generated structures, every such structure can be strongly embedded into $U$. Thus, U2 determines the Fra\"iss\'e limit uniquely.
\end{remark}


Now we consider a stronger amalgamation property known as the \emph{asymmetric amalgamation property}. However, in our examples the class $\mathfrak{C}_{f.g.}$ does not have this property, so we need to assume a subclass has that property. 

\begin{assumption}
Assume there is a subclass $\hat{\mathfrak{C}} \subseteq \mathfrak{C}$ with the following properties\footnote{$\hat{\mathfrak{C}}_{f.g.}$ denotes the collection of structures from $\hat{\mathfrak{C}}$ that are finitely $\hat{\mathfrak{C}}$-generated.}.
\begin{enumerate}
\item[{C3}] Every structure $A \in \mathfrak{C}$ has a unique (up to isomorphism over $A$) strong extension $\hat{A} \in \hat{\mathfrak{C}}$ which  is $\hat{\mathfrak{C}}$-generated by $A$. If $A \in \mathfrak{C}_{f.g.}$ then $\hat{A} \in \hat{\mathfrak{C}}_{f.g.}$.

\item[{C4}] If $A, B \in \mathfrak{C}$ with a strong embedding $A \hookrightarrow B$ then it can be extended to a strong embedding $\hat{A} \hookrightarrow \hat{B}$. 

\item[{C5}] $\hat{\mathfrak{C}}$ is closed under unions of countable strong chains.

\item[{AAP}] (Asymmetric Amalgamation Property) If $A_0,A_1,A_2 \in \hat{\mathfrak{C}}_{f.g.}$ with a strong embedding $A_0 \leq A_1$ and an embedding $A_0 \hookrightarrow A_2$ (not necessarily strong), then there is $B \in \hat{\mathfrak{C}}_{f.g.}$ with an embedding $A_1 \hookrightarrow B$ and a strong embedding $A_2 \leq B$ such that the corresponding diagram commutes. Moreover, if $A_0$ is strong in $A_2$ then $A_1$ is strong in $B$.
\end{enumerate}
\end{assumption}

\begin{proposition}\label{strong-amalgam}
If $\hat{\mathfrak{C}}$ satisfies \emph{AAP} then $\mathfrak{C}_{f.g.}$ has the strong amalgamation property.
\end{proposition}
\begin{proof}
Let $A, B_1, B_2 \in \mathfrak{C}_{f.g.}$ with strong embeddings $A \hookrightarrow B_1,~ A \hookrightarrow B_2$. By our assumptions we have strong extensions $A \leq \hat{A},~ B_1 \leq \hat{B}_1,~ B_2 \leq \hat{B}_2$ with $\hat{A},\hat{B}_1,\hat{B}_2 \in \hat{\mathfrak{C}}$. Moreover, there are strong embeddings of $\hat{A}$ into $\hat{B}_1$ and $\hat{B}_2$. Now we can use the AAP property\footnote{One actually does not need full AAP here. Strong amalgamation for $\hat{\mathfrak{C}}$ would be enough. Nevertheless, AAP will prove useful later on.} of $\hat{\mathfrak{C}}$ to construct a strong amalgam $B'$ of $\hat{B}_1$ and $\hat{B}_2$ over $\hat{A}$. Let $B$ be the substructure of $B'$ $\mathfrak{C}$-generated by $B_1$ and $B_2$. Clearly $B \in \mathfrak{C}_{f.g.}$ and it is a strong amalgam of $B_1$ and $B_2$ over $A$.
\end{proof}

\begin{notation}
For $A \in \hat{\mathfrak{C}}$ and a subset $X \subseteq A$, the substructure of $A$ $\mathfrak{C}$-generated by $X$ will be denoted by $\langle X \rangle^{\mathfrak{C}}_A$ while $\langle X \rangle^{\hat{\mathfrak{C}}}_A$ stands for the substructure of $A$ $\hat{\mathfrak{C}}$-generated by $X$. The same pertains to strong substructures generated by $X$ in the two classes. When no confusion can arise, we will drop the superscript.
\end{notation}

\begin{proposition}\label{same-strong-limit}
Under the assumptions \emph{C1-5, AAP}, the classes $\mathfrak{C}_{f.g.}$ and $\hat{\mathfrak{C}}_{f.g.}$ are strong amalgamation classes and have the same strong Fra\"iss\'e limit.
\end{proposition}
\begin{proof}
Firstly, we show that $\hat{\mathfrak{C}}$ is a strong amalgamation class. For this we need to prove that every countable $A \in \hat{\mathfrak{C}}$ has at most countably many strong finitely generated extensions in $\hat{\mathfrak{C}}$, up to isomorphism.

Let $B \in \hat{\mathfrak{C}}_{f.g.}$ be generated by $\bar{b}$ over $A$ as a $\hat{\mathfrak{C}}$-structure. Denote $B_0:=\lceil A \bar{b} \rceil^{\mathfrak{C}}_B$.  Then $B_0 \leq B$ and $B=\langle B_0 \rangle_B^{\hat{\mathfrak{C}}}$ which shows that $B=\hat{B_0}$. Since $A \leq B$, we have $A \leq B_0$ and so there are countably many choices for $B_0$ and hence countably many choices for $B$.

Let $U$ be the strong Fra\"iss\'e limit of $\hat{\mathfrak{C}}$. We will show that it satisfies U2 for $\mathfrak{C}_{f.g.}$. Let $A, B \in \mathfrak{C}_{f.g.}$ with strong embeddings $f:A \hookrightarrow B$ and $g:A \hookrightarrow U$. We can extend $f$ and $g$ to strong embeddings $\hat{A} \hookrightarrow \hat{B}$ and $\hat{A} \hookrightarrow U$ over $A$. Therefore $\hat{B}$ can be strongly embedded into $U$ over $\hat{A}$. The restriction of this embedding to $B$ will be a strong embedding of $B$ into $U$ over $A$.

Thus $U$ is also strongly saturated for $\mathfrak{C}$, hence $U$ is isomorphic to the Fra\"iss\'e limit of $\mathfrak{C}$.
\end{proof}


\begin{proposition}\label{ARP}
Under the above assumption $U$ has the following \emph{Asymmetric Richness Property}.
\begin{itemize}
\item[\emph{ARP}] If $A\leq B \in \hat{\mathfrak{C}}_{f.g.}$ then any embedding $A \hookrightarrow U$ extends to an embedding $B\hookrightarrow U$. Moreover, if the former embedding is strong then so is the latter.
\end{itemize}
\end{proposition}
\begin{proof}
Let $\lceil A \rceil \in \hat{\mathfrak{C}}_{f.g.}$ be the self-sufficient closure of $A$ in $U$ (in the sense of $\hat{\mathfrak{C}}$). By AAP there is $B' \in \hat{\mathfrak{C}}_{f.g.}$ with embeddings $\lceil A \rceil \leq B'$ and $B \hookrightarrow B'$ over $A$. Now richness of $U$ implies the desired result.
\end{proof}

ARP states that the amalgam $U$ is existentially closed in strong extensions, which is normally used to give a first-order axiomatisation of the amalgam.

In general U1 and U2 are not first-order axiomatisable, nor is ARP. normally they are $\mathfrak{L}_{\omega_1,\omega}$-axiomatisable provided the predimension has some definability properties (which we specify below). In order to extract a first-order axiomatisation from this $\mathfrak{L}_{\omega_1,\omega}$-axiomatisation, one normally approximates U1 and U2 by finitary axioms which are first-order. Wagner considers this problem in \cite{Wagner} and gives the appropriate conditions under which it can be done, working in a relational language though. In particular, if the language is finite and relational and $\mathfrak{C}_{f.g.}$ consists of finite structures then one can find a first-order axiomatisation of the amalgam. In general it is possible to give a similar first-order axiomatisation of $\Th(U)$ imposing quite strong definability conditions on $\delta$. However it seems those conditions would fail for the Ax-Schanuel predimension (see Section \ref{Predim-examples}) and so we consider weaker definability conditions.


Let $M \in \mathfrak{C}$ be an arbitrary structure.

\begin{definition}\label{definable-predim}
We say $\delta$ is \emph{(infinitely) definable} in $M$ if for any $n,m \in \mathbb{N}$ the set $\{ \bar{a}\in M^n: \delta(\bar{a})\geq m\}$ is definable by a possibly infinite Boolean combination of first-order formulas, i.e. an $\mathfrak{L}_{\omega_1,\omega}$-formula of the form
\begin{equation}
\bigwedge_{i<\omega}\bigvee_{j<\omega} \varphi_{i,j}^{m,n}(\bar{x}),
\end{equation}
where $\varphi_{i,j}^{m,n}(\bar{x})$ are first-order formulae. We say $\delta$ is \emph{universally definable} if the formulas $\varphi_{i,j}^{m,n}$ can be chosen to be universal formulas.
\end{definition}


Recall that we assumed $\delta$ is non-negative. This means, in particular, that
\begin{equation}\label{predim-ineq}
\delta(\bar{x})\geq 0 \mbox{ for all finite tuples } \bar{x} \subseteq M.
\end{equation}

\begin{lemma}
If $M \in \mathfrak{C}$ is saturated and $\delta$ is definable then the inequality \eqref{predim-ineq} is first-order axiomatisable.
\end{lemma}
\begin{proof}
By \eqref{predim-ineq} we know that for each $i$ we have
$$M \models \forall \bar{x} \bigvee_{j<\omega} \varphi_{i,j}^{0,n}(\bar{x}).$$
Since $M$ is saturated, there is a positive integer $N_i$ such that
$$M \models \forall \bar{x} \bigvee_{j<N_i} \varphi_{i,j}^{0,n}(\bar{x}).$$
Then \eqref{predim-ineq} is axiomatised by the following collection of axioms:
$$\forall \bar{x} \bigvee_{j<N_i} \varphi_{i,j}^{0,n}(\bar{x}),~ i<\omega, n<\omega.$$
\end{proof}

For a finite set $\bar{a} \subseteq M$ we say $\bar{a}$ is strong in $M$ if $\langle \bar{a} \rangle \leq M$. Definability of $\delta$ implies that for a finite set being strong in $M$ is $\mathfrak{L}_{\omega_1,\omega}$-definable. 


\begin{lemma}\label{near-model-complete}
Assume $U$ is saturated and $\delta$ is universally definable in $U$. Then $\Th(U)$ is nearly model complete, that is, every formula is equivalent to a Boolean combination of existential formulas in $U$.
\end{lemma}
\begin{proof}
For a finite tuple $\bar{a} \subseteq U$ its type (in $U$) is determined by the isomorphism type of $\lceil \bar{a} \rceil_U$ which is determined by finitely generated non-strong extensions of $\langle \bar{a} \rangle$ in $U$. If $\bar{a}$ and $\bar{b}$ satisfy exactly the same existential formulae (and hence exactly the same universal formulae), then for any non-strong extension of $\langle \bar{a} \rangle$ there is an isomorphic non-strong extension of $\langle \bar{b} \rangle$. Hence $\lceil \bar{a} \rceil_U \cong \lceil \bar{b} \rceil_U$. Thus, $\tp(\bar{a})$ is determined by existential formulae and their negations that are true of $\bar{a}$. Therefore $\Th(U)$ is nearly model complete.
\end{proof}

{When one knows the first-order theory of $U$, one can normally understand whether $U$ is saturated or not. It is saturated in our main examples, i.e. the exponential differential equation and the equation of the $j$-function (see Section \ref{j-chapter}). However, in general, it is possible to have a non-saturated Fra\"iss\'e limit. Baldwin and Holland \cite{Baldwin-Holland} give a criterion (called \emph{separation of quantifiers}) for saturatedness of $U$ (working under stronger definability conditions for $\delta$ though).}

\begin{definition}
We say $\delta$ is trivial on $\hat{\mathfrak{C}}$ if all embeddings of structures from $\hat{\mathfrak{C}}$ are strong.
\end{definition}

Note that in general $\delta$ is not defined on $\hat{\mathfrak{C}}$ (nor on $\hat{\mathfrak{C}}_{f.g.}$), so to be more precise we could say that strong embeddings induced by $\delta$ are trivial on $\hat{\mathfrak{C}}$. This is a weaker condition than triviality in the sense of Definition \ref{defin-predim-trivial-1}. \textbf{\emph{From now on, triviality of $\delta$ should be understood in this sense.}}

\begin{proposition}\label{predim-model-complete}
Assume $U$ is saturated. If $\delta$ is non-trivial on $\hat{\mathfrak{C}}$ then $\Th(U)$ is not model complete.
\end{proposition}
\begin{proof}
Non-triviality of the predimension means there are finitely generated $A \subseteq B \in \hat{\mathfrak{C}}_{f.g.}$ with $A\nleq B$. By universality of $U$ we know that there is a strong embedding of $A$ into $U$. Using the asymmetric amalgamation property we find a structure $U'\in \mathfrak{C}$ which extends $U$ and extends $B$ strongly such that the corresponding diagram commutes. This can be done since the amalgam $U$ is the union of a countable strong chain of finitely generated structures. So we can inductively use the asymmetric amalgamation for each of these structures and take the union of amalgams obtained in each step (these amalgams form a strong increasing chain). Then it is easy to see that $U \nleq U'$. On the other hand, $U'$ is countable and hence it can be embedded into $U$. Thus we have embeddings $U \hookrightarrow U' \hookrightarrow U$ and the first one is non-strong. Therefore we have a non-strong embedding of $U$ into itself. By Proposition \ref{elementary-strong} this embedding is not elementary which means $\Th(U)$ is not model complete.
\end{proof}

Now we define what it means for the inequality \eqref{predim-ineq} to be adequate.

\begin{definition}
Let $\hat{\mathfrak{C}} \subseteq  \mathfrak{C}$ be classes of structures closed under isomorphism and intersections and such that $\emptyset \in \mathfrak{C}$. Assume they satisfy FG, DCC, C1-5, AAP and $\delta$ is a non-negative universally definable predimension on $\mathfrak{C}_{f.g.}$. Let $M \in \mathfrak{C}$ be a countable structure. 
\begin{itemize}[leftmargin = 0.5cm]
\item We say that $\delta$ (or the inequality \eqref{predim-ineq}) is \emph{adequate} for $M$ if $U\equiv M$.
\item We say $\delta$ is \emph{strongly adequate} for $M$ if $M \cong U$.
\end{itemize}
\end{definition}

In other words, adequacy of a predimension inequality means that $\Th(M)$ can be obtained by a Hrushovski construction and strong adequacy means that the structure $M$ itself can be obtained by a Hrushovski construction. Strong adequacy also means that the reduct $\mathcal{K}_E$ is ``geometric'', that is, it is equipped with a pregeometry governed by the predimension. These notions will make more sense in differential setting where $M$ is always taken to be a reduct of a differentially closed field. Note also that when $M$ and $U$ are saturated, adequacy of $\delta$ implies its strong adequacy.

Note that we do not need definability of $\delta$ or AAP for some subclass $\hat{\mathfrak{C}}$ in order to construct the strong Fra\"iss\'e limit $U$ and define adequacy. However, these are natural assumptions since in most cases (in differential setting) the properties FG, C1-5 and definability of $\delta$ will be evident while strong amalgamation of $\mathfrak{C}_{f.g.}$ will be deduced from strong amalgamation of $\hat{\mathfrak{C}}_{f.g.}$, and in fact $\hat{\mathfrak{C}}_{f.g.}$ will have the asymmetric amalgamation property. That is the reason that we included all those conditions in the definition of adequacy. This will be illustrated in Section \ref{j-chapter}. 

\subsection{Examples}\label{Predim-examples}

In this section we give examples of predimensions that are the main motivating factor for this work.

\subsubsection{Complex exponentiation}

Let $\mathbb{C}_{\exp}:=(\mathbb{C};+,\cdot,0,1,\exp)$ be the complex exponential field. Let $E(x,y)$ be the graph of the exponential function and consider the structure $ \mathbb{C}_E:=(\mathbb{C};+,\cdot,0,1,E)$. Note that it is not saturated and its first-order theory is not stable since $\mathbb{Z}$ is definable.

For complex numbers $x_1, \ldots, x_n$ and their exponentials $y_1, \ldots, y_n$ define 
$$\delta(\bar{x},\bar{y}) := \td_\mathbb{Q}\mathbb{Q}(\bar{x}, \bar{y}) - \ldim_{\mathbb{Q}}(\bar{x}).$$ Schanuel's conjecture states non-negativity of this function.

Consider the class $\mathfrak{C}$ of all {(field-theoretically)} algebraically closed substructures of $\mathbb{C}_E$. For a finitely generated (i.e. of finite transcendence degree over $\mathbb{Q}$) substructure $A$ define
\begin{align*}
\sigma(A):=\max \{ n: &\mbox{ there are } a_i, b_i \in A,~ i=1,\ldots, n, \mbox{ with } a_i\mbox{'s}\\ & \mbox{ linearly independent over } \mathbb{Q} \mbox{ and } A \models E(a_i,b_i) \}
\end{align*}
and
$$\delta(A):=\td_{\mathbb{Q}}(A) - \sigma(A).$$

Then $\sigma$ is finite provided Schanuel's conjecture holds and $\delta$ is a well-defined non-negative predimension. 
However the inequality $\delta \geq 0$ is not first-order axiomatisable even assuming the conjecture holds. 

Schanuel's conjecture is wide open and so we cannot say much about this example. It is quite complicated from a model theoretic point of view. In particular, $\mathbb{Z}$ is definable in $\mathbb{C}_E$. So its first order theory is quite difficult to study. In spite of this Zilber discovered a nice way of treating the complex exponential field using infinitary logic. He considered algebraically closed fields with a relation which has some of the properties of complex exponentiation. Then he took all those structures where the analogue of Schanuel's conjecture holds. By a Hrushovski style construction he obtained a theory called \emph{pseudo-exponentiation}. It is axiomatised in the language $\mathfrak{L}_{\omega_1,\omega}(Q)$ where $Q$ is a quantifier for ``there exist uncountably many''. This theory (and its first-order part) is a natural candidate for the $\mathfrak{L}_{\omega_1,\omega}(Q)$-theory (respectively, first-order theory) of $\mathbb{C}_E$. Nevertheless, all these questions seem to be out of reach at the moment. We refer the reader to \cite{Zilb-pseudoexp,Kirby-Zilber-exp, Zilb-analytic-pseudoanalytic,Zilb-special-Schanuel,ZilbExpSUm, Kirby-pseudo} for details. Note also that many ideas in the analysis of the exponential differential equation (see below) originate in Zilber's work on pseudo-exponentiation. 

\begin{remark}
Submodularity does not hold for finite sets. Indeed, let $a, b \in \mathbb{C}$ with $\delta(a)=\delta(b)=1,~ \delta(a,b)=0$. Then taking $A=\{ a, b \}, B=\{ 2a, b\}$ we get
$$\delta(A \cup B) + \delta (A \cap B) = 0+1 > 0+0 = \delta(A)+\delta(B).$$
\end{remark}

\subsubsection{Exponential differential equation}

First, let us state the Ax-Schanuel theorem.

\begin{theorem}[Ax-Schanuel, \cite{Ax}]
Let $(K;+,\cdot,',0,1)$ be a differential field with field of constants $C$. If $(x_1,y_1),\ldots,(x_n,y_n)\in K^2$ are non-constant solutions to the exponential differential equation $y' = yx'$ then
\begin{equation}\label{Ax-Schanuel}
\td_CC(\bar{x},\bar{y})-\ldim_{\mathbb{Q}}(\bar{x}/C) \geq 1,
\end{equation}
where $\ldim_{\mathbb{Q}}(\bar{x}/C)$ is the dimension of the $\mathbb{Q}$-span of $x_1, \ldots, x_n$ in the quotient vector space $K/C$.
\end{theorem}  

Now let $\mathcal{K}:=(K;+,\cdot,',0,1)$ be a countable saturated differentially closed field with field of constants $C$. 
Let $\Exp(x,y)$ be defined by the exponential differential equation $y'=yx'$ and denote $\mathcal{K}_{\Exp}:=(K;+,\cdot,\Exp,0,1)$. 
Fix the language $\mathfrak{L}_{\Exp}:=\{ +, \cdot, \Exp, 0, 1\}$. Consider the following axioms for an $\mathfrak{L}_{\Exp}$-structure $F$ ($\mathbb{G}_a$ and $\mathbb{G}_m$ denote the additive and multiplicative groups of a field and $G_n:=\mathbb{G}_a^n \times \mathbb{G}_m^n$).

\begin{itemize}
\item[A1] $F$ is an algebraically closed field of characteristic $0$.

\item[A2] $C_F:= \{ c \in F: F \models \Exp(c,1) \}$ is a algebraically closed subfield of $F$.

\item[A3] $\Exp(F) = \{ (x,y)\in F^2: \Exp(x,y)\}$ is a subgroup of $G_1(F)$ containing $G_1(C_F)$.

\item[A4] The fibres of $\Exp$ in $\mathbb{G}_a(F)$ and $\mathbb{G}_m(F)$ are cosets of the subgroups $\mathbb{G}_a(C_F)$ and $\mathbb{G}_m(C_F)$ respectively.

\item[AS] For any $x_i, y_i \in F,~ i=1,\ldots,n$, if $F \models \bigwedge_{i=1}^n\Exp(x_i,y_i)$ and $\td_{C_F}(\bar{x},\bar{y}/C_F) \leq n$ then there are integers $m_1,\ldots,m_n$, not all of them zero, such that $m_1x_1+\ldots+m_nx_n \in C_F$.
\item[NT] $F \supsetneq C$.
\end{itemize}

Note that AS can be given by a first-order axiom scheme. A compactness argument gives a uniform version of AS. That is, given a parametric family of varieties $V(\bar{c})$ over $C$, there is a finite number $N$, such that if for some $\bar{c}$ we have $(\bar{x},\bar{y})\in V(\bar{c})$ and $\dim V(\bar{c})\leq n$ then $m_1x_1+\ldots+m_nx_n \in C$ for some integers $m_i$ with $|m_i| \leq N$.

Let $T^0_{\Exp}$ be the theory axiomatised by A1-A4, AS. The class $\mathfrak{C}$ consists of all countable models of $T^0_{\Exp}$ with a fixed field of constants $C$ (which is a countable algebraically closed field of transcendence degree $\aleph_0$). For $F \in \mathfrak{C}$ and $X \subseteq F$ we have $\langle X \rangle=C(X)^{\alg}$ with the induced structure from $F$. A structure $A \in \mathfrak{C}$ is finitely generated if and only if it has finite transcendence degree over $C$.

For finite tuples $\bar{x}, \bar{y} \in K^n$ with $\Exp(x_i,y_i)$ define 
$$\delta(\bar{x}, \bar{y}) := \td_{C}C(\bar{x},\bar{y}) - \ldim_{\mathbb{Q}}(\bar{x}/C).$$
The Ax-Schanuel theorem states positivity of this function (for non-constant tuples). It is easy to see that $\delta$ is universally definable.
We want to extend $\delta$ to $\mathfrak{C}_{f.g.}$. Following \cite{Kirby-semiab} for $A \in \mathfrak{C}_{f.g.}$ define 
\begin{align*}
\sigma(A):=\max \{ n: &\mbox{ there are } a_i, b_i \in A,~ i=1,\ldots, n, \mbox{ with } a_i\mbox{'s}\\ & \mbox{ linearly independent over } \mathbb{Q} \mbox{ mod } C \mbox{ and } A \models \Exp(a_i,b_i) \}
\end{align*}
and
$$\delta(A):=\td_{C}(A) - \sigma(A).$$

Firstly, note that $\sigma$ is well defined and finite since the Ax-Schanuel inequality bounds the number $n$ in consideration by $\td_{C}C(\bar{a},\bar{b})$ which, in its turn, is bounded by $\td_{C}A$.

Secondly, it is quite easy to prove that for any $A, B \in \mathfrak{C}_{f.g.}$
$$\sigma(A \cup B) \geq \sigma(A)+\sigma(B)-\sigma(A\cap B).$$
This implies that $\delta$ is submodular. Invariance of $\delta$ under isomorphism is clear too. Hence it is a predimension.

The Ax-Schanuel inequality is equivalent to saying that $\delta(A) \geq 0$ for all $A \in \mathfrak{C}_{f.g.}$ where equality holds if and only if $A=C$.

The class $\mathfrak{C}$ satisfies the strong amalgamation property but not the asymmetric amalgamation property. So we let $\hat{\mathfrak{C}}$ be the subclass of $\mathfrak{C}$ consisting of \emph{full} structures. A structure $A \in \mathfrak{C}$ is \emph{full} if for every $a \in A$ there are $b_1,b_2 \in A$ with $A \models \Exp(a,b_1) \wedge \Exp(b_2,a)$. Then $\hat{\mathfrak{C}}$ has the AAP property and satisfies all the assumptions made in the previous sections.


\begin{theorem}[\cite{Kirby-semiab}]\label{Kirby-adequate}
The Ax-Schanuel inequality is strongly adequate for $\mathcal{K}_{\Exp}$. 
\end{theorem}


Let us give a complete axiomatisation of $\Th(\mathcal{K}_{\Exp})$. For that we will need to formulate an existential closedness statement. For a $k \times n$ matrix $M$ of integers we define $[M]:G_n(F) \rightarrow G_k(F)$ to be the map given by $[M]:(\bar{x},\bar{y}) \mapsto (u_1,\ldots,u_k, v_1,\ldots, v_k)$ where
$$u_i = \sum_{j=1}^n m_{ij}x_j \mbox{ and } v_i = \prod_{j=1}^n y_j^{m_{ij}}.$$

\begin{definition}\label{rotund}
An irreducible variety $V \subseteq G_n(F)$ is \emph{rotund} if for any $1 \leq k \leq n$ and any $k\times n$ matrix $M$ of integers $\dim [M](V) \geq \rank M$. If for any non-zero $M$ the stronger inequality $\dim [M](V) \geq \rank M + 1$ holds then we say $V$ is \emph{strongly rotund}.
\end{definition}

The definition of rotundity is originally due to Zilber though he initially used the word \textit{normal} for these varieties \cite{Zilb-pseudoexp}. The term \emph{rotund} was coined by Kirby in \cite{Kirby-semiab}. 

Strong rotundity fits with the Ax-Schanuel inequality in the sense that it is a sufficient condition for a variety defined over $C$ to contain a non-constant exponential point. More precisely, if $\mathcal{F}$ is differentially closed and $V \subseteq G_n(F)$ is a strongly rotund variety defined over the constants, then the intersection $V(F) \cap \Exp(F)$ contains a non-constant point. 

Nevertheless, the existential closedness axiom we will use for the axiomatisation of $T_{\Exp}$ is slightly different. One needs to consider varieties that are not necessarily defined over $C$.

The existential closedness property for a model $F$ of $T_{\Exp}^0$ is as follows.

\begin{itemize}
\item[EC] For each irreducible rotund variety $V \subseteq G_n(F)$ the intersection $V(F) \cap \Exp(F)$ is non-empty.
\end{itemize}
As noted above, $V$ is not necessarily defined over $C$ and the point in the intersection may be constant. 


Rotundity of a variety is a definable property. This allows one to axiomatise the above statement by a first-order axiom scheme. Reducts of differentially closed fields satisfy EC and it gives a complete theory together with the axioms mentioned above.

\begin{theorem}[\cite{Kirby-semiab}]\label{exp-axioms}
The first-order theory of an exponential reduct of a differentially closed field is axiomatised by the following axioms and axiom schemes: \emph{A1-A4, AS, EC, NT}.
\end{theorem}

In \cite{Aslanyan-linear} we generalised the Ax-Schanuel theorem to linear differential equations of arbitrary order with constant coefficients and established the adequacy of those predimension inequalities.

In Section \ref{j-chapter} we study the predimension given by the Ax-Schanuel inequality for the $j$-function and give full details of the construction and axiomatisation of the Fra\"{i}ss\'{e} limit.

\subsection{Predimensions in the differential setting}\label{dif-eq}

Let $\mathcal{K}:=(K;+,\cdot,',0,1)$ be a countable saturated differentially closed field with field of constants $C$. 
Suppose $f(X,Y) \in \mathbb{Q}\{ X, Y\}$ is a differential polynomial with $\ord_Y(f)=m+1$. Consider the differential equation
\begin{equation}\label{equation}
f(x,y)=0.
\end{equation}

Let $E(x,y_0,\ldots,y_m)$ be an $(m+2)$-ary relation defined by 
$$ f(x,y_0)=0 \wedge \bigwedge_{i=0}^{m-1}  y'_{i}=y_{i+1}x'.$$

If $x$ is non-constant and $E(x,y_0,\ldots,y_m)$ holds then $y_i = \partial_x^i y_0$ where $\partial_x$ is the differentiation with respect to $x$ given by $u \mapsto \frac{u'}{x'}$. Thus we think of $y_i$ as the $i$-th derivative of $y_0$ with respect to $x$.

We fix the language $\mathfrak{L}_E:=\{ +, \cdot, E, 0, 1 \}$. Let $\mathfrak{C}$ be a class of $\mathfrak{L}_E$-structures satisfying all requirements set in Sections \ref{predimensions} and \ref{amalgamation} (in particular, the existence of $\hat{\mathfrak{C}}$ with the appropriate properties is assumed). Assume $\delta$ is a non-negative predimension on $\mathfrak{C}_{f.g.}$. normally $\mathfrak{C}$ will consist of algebraically closed fields with a relation $E$ satisfying some basic universal axioms of $E$-reducts of differential fields. These axioms will depend on functional equations satisfied by $E$. Most importantly, we should have an axiom scheme for the inequality $\delta \geq 0$. 


\begin{definition}
We say $\delta$ is (strongly) adequate (for the differential equation $E$) if it is (strongly) adequate for the reduct $\mathcal{K}_E:=(K;+,\cdot,E,0,1)$.
\end{definition}

\begin{remark}
It makes sense to consider just a binary relation for the set of solutions of our differential equation, without including derivatives, and study predimensions in that setting. More generally, we can do the same for an arbitrary reduct of a differentially closed field and define adequacy as above.
\end{remark}


Now we consider a special kind of predimension {motivated by the Ax-Schanuel inequality for the exponential differential equation and its analogue for the $j$-function}. Assume $d$ is a modular dimension function on $\mathcal{K}$. Suppose whenever $(x_i,y_i)$ are solutions of equation \eqref{equation}, the following inequality holds: {
\begin{equation}\label{inequality}
\td_CC( \{ x_i, \partial^j_{x_i}y_i: i=1,\ldots,n,~ j=0,\ldots,m\}) - (m+1)d(\bar{x},\bar{y})\geq 0.
\end{equation} }

The inequality \eqref{inequality} is first-order axiomatisable provided that $d$ is type-definable in the algebraically closed field $K$, i.e. for each $m$ and $n$ the set $\{ \bar{x} \in K^n: d(\bar{x}) \geq m \}$ is type definable (in the language of rings).


For $A \in \mathfrak{C}_{f.g.}$ define 
\begin{align*}
\sigma(A):=\max \{ & d(\bar{a},\bar{b}): a_i, b_i \in A \mbox{ and there are }  b^1_i,\ldots, b^m_i \in A,\\ &\mbox{with } A \models E(a_i,b_i,b^1_i, \ldots, b^m_i) \}
\end{align*}
and
$$\delta(A):=\td(A/C) - (m+1) \cdot \sigma(A).$$


It is easy to see that $\sigma$ is finite and hence $\delta$ is well defined. On the other hand for $A, B \in \mathfrak{C}_{f.g.}$ one can easily prove (using modularity of $d$) that
$$\sigma(A \cup B) \geq \sigma(A)+\sigma(B)-\sigma(A\cap B).$$
Thus, $\delta$ is submodular. In this manner we obtain a predimension on $\mathfrak{C}_{f.g.}$ and it makes sense to ask whether it is adequate or not. 

As pointed out above, adequacy gives a good understanding of the structure of the appropriate reduct of a differentially closed field. Let us elaborate on this intuitive idea, based on the analysis of pseudo-exponentiation and the exponential differential equation (and the differential equation of the $j$-function in Section \ref{j-chapter}).

In order to understand the structure of our differential equation, one has to understand which systems of equations in the language of the reduct $\mathcal{K}_E$ have a solution. Then a predimension inequality (like \eqref{inequality}) implies that ``overdetermined'' systems cannot have solutions. Adequacy means that this is the only obstacle: if having a solution does not contradict our inequality then there is a solution. It is not difficult to see that this question is equivalent to understanding which varieties contain (generic enough) points that are solutions to our differential equation\footnote{We do not state precisely what we mean by this because we will see it in the case of the exponential differential equation and the equation of $j$ which will be enough to understand the question in general.} (we call them $E$-points). This is in fact how one axiomatises the first-order theory of a differential equation (i.e. the theory of the corresponding reduct) with an adequate predimension inequality. 

{Indeed, as we noted in Section \ref{amalgamation}, one normally approximates the richness property (which determines the strong Fra\"{i}ss\'{e} limit uniquely up to isomorphism) by first-order axioms in order to give an axiomatisation of $\Th(U)$. Richness of the strong Fra\"iss\'e limit $U$ implies that it is existentially closed in strong extensions. So if a variety contains an $E$-point in a strong extension of $U$ then such a point already exists in $U$. When one tries to axiomatise this property, one normally proves that varieties with certain properties always contain an $E$-point. However, according to the richness property, we need also make sure that when we work over a strong substructure as a set of parameters then there exists an $E$-point in our variety which is strongly embedded into $U$. So, our axioms should state that varieties with the appropriate properties contain an $E$-point which cannot be extended to another point with lower predimension.  In this case the axiomatisation is $\forall \exists \forall$.} 

{However, in our main examples, that is, the exponential differential equation and the equation of the $j$-function, we end up with simpler axioms which are in fact $\forall \exists$. Let us explain how one obtains those axioms. Suppose we work over a strong substructure $A \leq U$ and $V$ is a variety defined over $A$. If we know that $V$ contains an $E$-point $\bar{b}$ and $\delta(\bar{b}/A)>0$ then it is possible that $\bar{b}$ is not strong in $U$. This can happen if $V$ has high dimension. In such a situation one uses the tool of intersecting varieties with generic hyperplanes (see Lemma \ref{generic-hyperplanes}) and decreases the dimension of $V$, more precisely, one replaces $V$ with a subvariety $V'$ defined over some $A'$ with $A \leq A' \leq U$. Now if $\dim V'$ is small enough then an $E$-point $\bar{b}$ in $V'$ satisfies $\delta(\bar{b}/A') =0$ which shows that $\bar{b}$ is strong in $U$ (since $A'\leq U$). Thus, the existence of $E$-points in certain varieties is enough to deduce the existence of $E$-points which are strong in $U$. Hence one axiomatises the existential closedness property by saying that certain varieties contain $E$-points. Then one normally ends up with an $\forall \exists$ axiom scheme which, along with the basic universal axioms (including an axiom scheme stating non-negativity of $\delta$), is expected to give a complete axiomatisation of the theory of the strong Fra\"iss\'e limit. So, in this case the axiomatisation is expected to be $\forall \exists$.} 

{This observation justifies the condition of $\forall \exists$-axiomatisability in Theorem \ref{inductive}. Nevertheless, we recall once more that those speculations are based on the aforementioned examples, and in general we expect an $\forall \exists \forall$-axiomatisation rather than just $\forall \exists$. On the other hand, the procedure described above and in particular the method of intersecting varieties with generic hyperplanes is quite general and can be carried out for various differential equations with a predimension inequality. So in ``nice'' examples we hope to get an $\forall \exists$ theory. In Section \ref{j-chapter} we illustrate those ideas on the example of the differential equation of the $j$-function. Finally, let us remark that getting a first order axiomatisation for the Fra\"iss\'e limit is by no means ``automatic'' since some technical issues may arise depending on the setting as we will see in Section \ref{j-general-case}.}  

\section{Connection to definability of derivations in reducts of differentially closed fields}\label{connection-def-der}

Let $\mathcal{F}=(F;+,\cdot,0,1,D)$ be a differentially closed field. In \cite{Aslanyan-def-deriv} we considered the question of definability of the derivation $D$ in reducts of $\mathcal{F}$ of the form $\mathcal{F}_{R}=(F;+,\cdot,0,1,P)_{P \in R}$ where $R$ is some collection of definable sets in $\mathcal{F}$. It turns out that this question is closely related to the existence of an adequate predimension inequality in the appropriate reduct. Intuitively, if $D$ is definable in a reduct then finding an adequate predimension inequality in that reduct would mean that we can find an adequate predimension in a differentially closed field. However, such a predimension must be trivial. We give two precise results below that support this idea.

\begin{theorem}[\cite{Aslanyan-def-deriv}]\label{inductive}
If $T_{R}:=\Th(\mathcal{F}_R)$ is inductive (i.e. $\forall \exists$-axiomatisable) and defines $D$ then it is model complete.
\end{theorem}

\begin{corollary}
If $D$ is definable in a reduct $\mathcal{F}_E$ of $\mathcal{F}$ and $\Th(\mathcal{F}_E)$ has an $\forall \exists$-axiomatisation then $E$ cannot have a non-trivial strongly adequate predimension inequality.
\end{corollary}
\begin{proof}
This follows immediately from Theorem \ref{inductive} and Proposition \ref{predim-model-complete}
\end{proof}

\begin{corollary}
The exponential differential equation does not define $D$.
\end{corollary}


We have one more result in this direction. Let $\mathcal{F}$ be a countable saturated differentially closed field. Assume $\mathfrak{C}$ is a collection of structures in the language of reducts $\mathfrak{L}_R$ and $\delta$ is a predimension on $\mathfrak{C}_{f.g.}$ satisfying all necessary conditions given in Section \ref{predimensions-and-Hrush}. Let $d$ be the dimension associated with $\delta$. {Below by a $d$-generic type (over some parameter set $A$) we mean the type of an element $a$ with $d(a/A)=1$.}


\begin{theorem}\label{no-predim}
Assume the underlying fields of structures from $\mathfrak{C}_{f.g.}$ are algebraically closed of finite transcendence degree over $\mathbb{Q}$. Assume further that $d$-generic $1$-types (over finite sets) are not algebraic. If $\D$ is definable in $\mathcal{F}_R$ and $\delta$ is strongly adequate, then the reduct is model complete and hence $\delta$ is trivial.
\end{theorem}
In general, it is possible that a $d$-generic $1$-type is not unique. Moreover, in some trivial examples such a type may be algebraic. So our assumption excludes such degenerate cases. {In particular, if the free amalgamation property holds for $\mathfrak{C}_{f.g.}$ then $d$-generic types cannot be algebraic.} Actually, it will suffice to assume that generic $1$-types have more than one realisation. In fact, we expect $d$-generic types to be generic in the sense of the reduct of a differentially closed field. As we know those are unique and have maximal rank.

We will need the following lemma in the proof of the above theorem.

\begin{lemma}[\cite{Aslanyan-def-deriv}, Propositions 6.1 and 4.5]\label{existential-def}
Let $a\in F$ be a differentially transcendental element over $\mathbb{Q}$. Suppose $\varphi(x,y)$ is a formula in the language of reducts $\mathfrak{L}_R$ such that 
$$\mathcal{F} \models \forall y ( \varphi(a,y) \leftrightarrow y=D a).$$ 
Then $D$ is definable (without parameters). Moreover, if $\varphi$ is existential then $D$ is existentially definable and $\Th(\mathcal{F}_R)$.
\end{lemma}

\begin{proof}[Proof of Theorem \ref{no-predim}]
Strong adequacy means that $\mathcal{F}_R$ is the Fra\"{i}ss\'{e} limit of $\mathfrak{C}_{f.g.}$. Let $a \in F$ be differentially transcendental. 
Denote $A:= \lceil a \rceil$ (the strong closure of $a$ in $\mathcal{F}_R$) and $A':=\lceil a, \D a \rceil$. 

If $d(\D a/a)=1$ then by our assumption $\tp_R(\D a/a)$ (which is a $d$-generic type) has more than one realisation which contradicts definability of $\D a$ over $a$. Thus, $d(\D a/a)=0$ and so $\delta(A')=d(a, \D a) = d(a)=\delta(A)$. {Since $A \subseteq A'$, we have $\delta(A'/A)= \delta(A') - \delta(A) = 0$.} Let $A'$ be $\mathfrak{C}$-generated by $(a,\D a, \bar{u})$. Extending $\bar{u}$ if necessary we can assume that $A'=\mathbb{Q}(a,\D a,\bar{u})^{\alg}$ (here we use the fact that $\td(A'/\mathbb{Q}) < \aleph_0$). 

If $(v, \bar{w})$ is a realisation of the existential type $\etp_R(\D a, \bar{u}/a)$ then we claim that $v = \D a$. Indeed, in a differentially closed field the type of a (field-theoretically) algebraic element is isolated by its minimal algebraic equation (and so all algebraic conjugates of that element have the same type), hence 
$B:=\mathbb{Q}(a, v, \bar{w})^{\alg} \cong \mathbb{Q}(a, \D a, \bar{u})^{\alg}=A'$ where the isomorphism is in the sense of $\mathfrak{L}_R$-structures (induced from $\mathcal{F}_R$). Therefore $B \in \mathfrak{C}$ and $\delta(B)=\delta(A')$. If $d(a)=0$ then $\delta(A')=0$ and so $\delta(B)=0$. If $d(a)=1$ then $A=\langle a \rangle$ is the structure $\mathfrak{C}$-generated by $a$. Since $a \in B$, we must have $A \subseteq B$ and $\delta(B/A)=\delta(B)-\delta(A)=\delta(A')-\delta(A)=0$. In both cases $B \leq \mathcal{F}_R$. Now by homogeneity of the Fra\"{i}ss\'{e} limit for strong substructures, the above isomorphism between $B$ and $A'$ extends to an automorphism of $\mathcal{F}_R$. This implies $\tp_R(\D a, \bar{u}/a) = \tp_R(v, \bar{w}/a)$ and so $\tp_R(\D a/a) = \tp_R(v/a)$. On the other hand, $\D a$ is definable over $a$, hence we must have $v=\D a$.

Thus, for $p(x,y,\bar{z}):=\etp_R(a,\D a, \bar{u})$ we have 
$$\mathcal{F}_R \models \forall y \left( \exists \bar{z} \bigwedge p(a,y,\bar{z}) \longleftrightarrow y= \D a \right).$$

A compactness argument shows that there is an (existential) $\varphi(x,y,\bar{z}) \in p$ so that $$\mathcal{F}_R \models \forall y \left( \exists \bar{z} \varphi(a,y,\bar{z}) \longleftrightarrow y= \D a \right).$$ 
By Lemma \ref{existential-def}, $\D$ is existentially definable in $\mathcal{F}_R$. Therefore, by \cite[Proposition 4.5]{Aslanyan-def-deriv}, $\mathcal{F}_R$ is model complete.
\end{proof}

The result will still hold if instead of assuming that finitely generated structures have finite transcendence degree we assume $\delta$ is quantifier-free (infinitely) definable.

\section{The $j$-function}\label{j-chapter}

In this section we will study the Ax-Schanuel inequality for the $j$-function established by Pila and Tsimerman. We show that the models of a theory (which is essentially the universal theory of appropriate reducts of differential fields) have the strong amalgamation property (along with all other necessary properties), construct the strong Fra\"iss\'e limit $U$ and give an axiomatisation of its first-order theory. This gives an axiomatisation of the theory of the differential equation of the $j$-function assuming adequacy of the Ax-Schanuel inequality. We do not prove adequacy in this paper, but it follows from the results of this section and those of \cite{Aslanyan-Eterovic-Kirby-Diff-EC-j}.


The definitions and results of this section are analogous to their exponential counterparts. Many proofs are adapted from \cite{Kirby-semiab} and \cite{Bays-Kirby-exp}. However, we should note that some things are simpler for $j$ while others are subtler and more complicated. Moreover, there are some inaccuracies in \cite{Kirby-semiab} some of which have been corrected in various papers, but some still remain. For instance, the exponential analogue of Proposition \ref{ID} of this paper is missing from \cite{Kirby-semiab} (it is assumed to be trivial in the proof of \cite[Theorem 4.11]{Kirby-semiab}). Those issues are addressed in this section for the $j$-function, and their exponential counterparts can be dealt with similarly. 

\subsection{Background on the $j$-function}\label{jBackground}

We do not need to know much about the $j$-function itself, nor need we know its precise definition. Being familiar with some basic properties of $j$ will be enough for this section. We summarise those properties below referring the reader to \cite{Lang-elliptic,Serre,Masser-Heights,Silverman} for details.

Let $\GL_2(\mathbb{C})$ be the group of $2 \times 2$ matrices with non-zero determinant. This group acts on the complex plane (more precisely, the Riemann sphere) by linear fractional transformations. Namely, for a matrix $g=\begin{pmatrix}
a & b \\
c & d
\end{pmatrix} \in \GL_2(\mathbb{C})$ we define 
$$gz=\frac{az+b}{cz+d}.$$
The subgroup $\GL_2^+(\mathbb{R})$ of real matrices with positive determinant acts on the upper half-plane $\mathbb{H}:=\{ z \in \mathbb{C}: \im(z)>0 \}$. This induces similar actions of the subgroups of $\GL_2^+(\mathbb{R})$ such as $\SL_2(\mathbb{Z})$ and $\GL_2^+(\mathbb{\mathbb{Q}})$ which play an important role in the theory of modular forms.

The function $j$ is a modular function of weight $0$ for the \emph{modular group} $\SL_2(\mathbb{Z})$, which is defined and analytic on the upper half-plane $\mathbb{H}$. It is $\SL_2(\mathbb{Z})$-invariant, and by means of $j$ the quotient $Y(1)=\SL_2(\mathbb{Z})\setminus \mathbb{H}$ is identified with $\mathbb{C}$. Thus, $j$ is a bijection from the fundamental domain of $\SL_2(\mathbb{Z})$ to $\mathbb{C}$. 

For $g \in \GL_2^+(\mathbb{Q})$ we let $N(g)$ be the determinant of $g$ scaled so that it has relatively prime integral entries. For each positive integer $N$ there is an irreducible polynomial $\Phi_N(X,Y)\in \mathbb{Z}[X,Y]$ such that whenever $g \in \GL_2^+(\mathbb{Q})$ with $N=N(g)$, the function $\Phi_N(j(z),j(gz))$ is identically zero. Conversely, if $\Phi_N(j(x),j(y))=0$ for some $x,y \in \mathbb{H}$ then $y=gx$ for some $g \in \GL_2^+(\mathbb{Q})$ with $N=N(g)$. The polynomials $\Phi_N$ are called \emph{modular polynomials}. It is well known that $\Phi_1(X,Y)=X-Y$ and all the other modular polynomials are symmetric. For $w=j(z)$ the image of the $\GL_2^+(\mathbb{Q})$-orbit of $z$ under $j$ is called the \emph{Hecke orbit} of $w$. It obviously consists of the union of solutions of the equations $\Phi_N(w,X)=0,~ N\geq 1$. Two elements $w_1,w_2 \in \mathbb{C}$ are called \emph{modularly independent} if they have different Hecke orbits, i.e. do not satisfy any modular relation $\Phi_N(w_1,w_2)=0$. This definition makes sense for arbitrary fields (of characteristic zero) as the modular polynomials have integer coefficients.

The $j$-function satisfies an order $3$ algebraic differential equation over $\mathbb{Q}$, and none of lower order (i.e. its differential rank over $\mathbb{C}$ is $3$). Namely, $F(j,j',j'',j''')=0$ where 
$$F(y_0,y_1,y_2,y_3)=\frac{y_3}{y_1}-\frac{3}{2}\left( \frac{y_2}{y_1} \right)^2 + \frac{y_0^2-1968y_0+2654208}{2y_0^2(y_0-1728)^2}\cdot y_1^2.$$

Thus
$$F(y,y',y'',y''')=S(y)+R(y)(y')^2,$$
where $S$ denotes the \emph{Schwarzian derivative} defined by $S(y) = \frac{y'''}{y'} - \frac{3}{2} \left( \frac{y''}{y'} \right) ^2$ and $R(y)=\frac{y^2-1968y+2654208}{2y^2(y-1728)^2}$.

The following result is well known. The proof is taken from  \cite{Freitag-Scanlon}.
\begin{lemma}\label{solutions-of-eq-j}
All functions $j(gz)$ with $g \in \SL_2(\mathbb{R})$ satisfy the differential equation $$F(y,y',y'',y''')=0$$ and all solutions (defined on $\mathbb{H}$) are of that form. If we allow functions not necessarily defined on $\mathbb{H}$, then all solutions will be of the form $j(gz)$ where $g \in \SL_2(\mathbb{C})$.
\end{lemma}
\begin{proof}
Let $f(z)$ be a meromorphic function defined on some open domain $U\subseteq \mathbb{C}$. Since $j:\mathbb{H}\rightarrow \mathbb{C}$ is surjective, there is a function $h:U\rightarrow \mathbb{H}$ such that $j(h(z))=f(z)$ on $U$. Applying the Schwarzian derivative to this equality we get
$$S(f) = S(j \circ h) = (S(j)\circ h)\cdot (h')^2 + S(h).$$ Therefore $S(f) +R(f) \cdot (f')^2=0$ if and only if
$$ (S(j)\circ h)\cdot (h')^2 + S(h) + R(j\circ h)\cdot (j'\circ h)^2 \cdot (h')^2=0.$$
On the other hand we have
$S(j) +R(j) \cdot (j')^2=0$, hence $S(j)\circ h +R(j\circ h) \cdot (j'\circ h)^2=0$. Thus, $F(f,f',f'',f''')=0$ if and only if $S(h)=0$, i.e. $h=gz$ for some $g\in \SL_2(\mathbb{C})$.
\end{proof}

\subsection{Ax-Schanuel and Weak Modular Zilber-Pink}\label{jAS-ZP}

\begin{theorem}[Ax-Schanuel for $j$, \cite{Pila-Tsim-Ax-j}]\label{j-chapter-Ax-for-j}
Let $(K;+,\cdot,',0,1)$ be a differential field and let $z_i, j_i \in K \setminus C,~ j_i^{(1)}, j_i^{(2)}, j_i^{(3)} \in K^{\times},~ i=1,\ldots,n,$ be such that
\begin{equation*}
F\left(j_i,j_i^{(1)}, j_i^{(2)}, j_i^{(3)}\right)=0 \wedge j_i'=j_i^{(1)}z_i' \wedge \left(j_i^{(1)}\right)'=j_i^{(2)}z_i'\wedge \left(j_i^{(2)}\right)'=j_i^{(3)}z_i'.
\end{equation*}
If the $j_i$'s are pairwise modularly independent then 
\begin{equation}\label{j-chapter-Ax-ineq}
\td_CC\left(\bar{z},\bar{j},\bar{j}^{(1)},\bar{j}^{(2)}\right) \geq 3n+1.
\end{equation}
\end{theorem}

\begin{corollary}[Ax-Schanuel without derivatives]\label{j-Ax-without-der}
If $z_i, j_i$ are non-constant elements in a differential field $K$ with $F(j_i,\partial_{z_i}j_i,\partial^2_{z_i}j_i,\partial^3_{z_i}j_i)=0$ for $i=1,\ldots,n$, then 
$$\td_CC(\bar{z},\bar{j}) \geq n+1,$$
unless for some $N, i, k$ we have $\Phi_N(j_i,j_k)=0$.
\end{corollary}

This theorem implies in particular that the only algebraic relation between the functions $j(z)$ and $j(gz)$ for $g \in \SL_2^+(\mathbb{R})$ are the modular relations (corresponding to $g \in \GL^+_2(\mathbb{Q})$).

\iftrue
An important consequence of the Ax-Schanuel theorem is a weak form of the modular Zilber-Pink conjecture. Below $K$ is an algebraically closed field.

\begin{definition}\label{def-special}
A \emph{special variety} is an irreducible component of a Zariski-closed set defined by some modular equations. Note that we allow a modular equation of the form $\Phi_N(x_i,x_i)=0$ which is equivalent to allowing equations of the form $x_i=c$ where $c$ is a special point (the image of a quadratic number under $j$).
\end{definition}

\begin{definition}
Let $V \subseteq K^n$ be an algebraic variety. An \emph{atypical subvariety} of $V$ is an irreducible component $W$ of some $V \cap S$, where $S$ is a special subvariety, such that $$\dim W> \dim V + \dim S -n.$$
An atypical subvariety $W$ of $V$ is said to be \emph{strongly atypical} if it is not contained in any hyperplane of the form $x_i=a$ for some $a \in K$ (i.e. no coordinate is constant on $W$).
\end{definition}

Note that in general when $U$ is a smooth variety and $V, W \subseteq U$ are subvarieties then
$$\dim (V \cap W) \geq \dim V + \dim W - \dim U.$$
For varieties in ``general position'' we expect that equality should hold. However, when the defining algebraic equations of $V$ and $W$ are not ``independent'' (inside $U$) then we may have a strict inequality. In other words, the intersection is atypically large. This is the motivation behind the above definition.

The following is an analogue of Zilber's Conjecture on Intersections with Tori (see \cite{Zilb-exp-sum-published,Kirby-Zilber-exp}).

\begin{conjecture}[Modular Zilber-Pink; \cite{Pila-Tsim-Ax-j}]\label{modular-ZP}
Every algebraic variety contains only finitely many maximal atypical subvarieties.
\end{conjecture}

\begin{definition}
When $V \subseteq K^{n+m}$ is a variety defined over $\mathbb{Q}$, and $A \subseteq K^m$ is its projection onto the last $m$ coordinates ($A$ is a constructible set), for each $\bar{a} \in A$ we let $V_{\bar{a}}$ (or $V(\bar{a})$) be the fibre of the projection above $\bar{a}$. The family $(V_{\bar{a}})_{\bar{a} \in A}$ is called a \emph{parametric family of varieties}.
\end{definition}

The following theorem is a weak version of the modular Zilber-Pink conjecture and follows from Ax-Schanuel. Pila and Tsimerman \cite{Pila-Tsim-Ax-j} give an o-minimality proof. In \cite{Aslanyan-weakMZPD} we give a differential algebraic proof adapting the proof of weak CIT by Zilber \cite{Zilb-exp-sum-published} and Kirby \cite{Kirby-semiab}.

\begin{theorem}[Weak modular Zilber-Pink]\label{weak-modularZP}
Let $K$ be an algebraically closed field. Given a parametric family of algebraic varieties $(V_{\bar{a}})_{\bar{a}\in A}$ in $K^n$, there is a finite collection of proper special varieties $(S_i)_{i\leq N }$ in $K^n$ such that for every $\bar{a} \in A$, every strongly atypical subvariety of $V_{\bar{a}}$ is contained in one of $S_i$.
\end{theorem}
\fi

\subsection{The universal theory and predimension}\label{j-univtheory-predim}
For simplicity, we are going to work mainly with Ax-Schanuel without derivatives. However, most of our results remain true if we consider derivatives too, and in the last section we will formulate definitions and main results in that generality, pointing out an issue related to weak modular Zilber-Pink ``with derivatives''.


Let $(K;+,\cdot,',0,1)$ be a differential field and let $F(y,y',y'',y''')=0$ be the differential equation of the $j$-function. 
Consider its two-variable version\footnote{Recall that for a non-constant $x$ we define $\partial_x:y \mapsto \frac{y'}{x'}$.}
\begin{equation}
f(x,y):=F\left(y,\partial_xy,\partial^2_xy,\partial^3_xy\right)=0.
\end{equation}



We prove several lemmas about this differential equation. Below $C$ denotes the field of constants of $K$.

\begin{lemma}\label{j-z-fibre}
Given $a_1,a_2,b \in K \setminus C$, if $f(a_1,b)=0$ then $f(a_2,b)=0$ iff $a_2=ga_1$ for some $g \in \SL_2(C)$.
\end{lemma}
\begin{proof}
By replacing the field derivation with $\partial_{a_1}$, we may assume without loss of generality that $a_1=t$ with $t'=1$. For simplicity write $a_2=:a$. Let also $S_{\partial_a}$ be the Schwarzian derivative with respect to $\partial_a$. Then we know that $Sb+R(b)(b')^2=0$, and so $S_{\partial_a}(b)+R(b)(\partial_ab)^2=0$ if and only if $(a')^2\cdot S_{\partial_a}(b) = S(b)$. However, straightforward calculations show that $(a')^2\cdot S_{\partial_a}(b) = S(b) - S(a)$. Hence, $f(a,b)=0$ iff $S(a)=0$ iff $a = gt$ for some $g \in \SL_2(C)$.
\end{proof}

\begin{lemma}
If $f(z,j_1)=0$ for some non-constants $z, j_1$, and $j_2$ satisfies $\Phi_N(j_1,j_2)=0$ for some modular polynomial $\Phi_N$ then $f(z,j_2)=0$.
\end{lemma}
\begin{proof}
Let $K_0 := \mathbb{Q} \langle z, j_1, j_2 \rangle$ be the differential subfield of $K$ generated by $z, j_1, j_2$, and let $C_0$ be its field of constants. Embedding $K_0$ into the field of germs of meromorphic functions, we can assume $C_0 \subseteq \mathbb{C}$ and $j_1, j_2$ are complex meromorphic functions of variable $z$. But then {by Lemma \ref{solutions-of-eq-j}}, $j_1=j(g_1z)$ for some $g_1 \in \SL_2(\mathbb{C})$ where $j: \mathbb{H} \rightarrow \mathbb{C}$ is the $j$-invariant. Now the identity $\Phi_N(j_1(z),j_2(z))=0$ implies $j_2(z)=j(g_2z)$ where $g_2=gg_1$ for some $g \in \GL_2^+(\mathbb{Q})$. {Applying Lemma \ref{solutions-of-eq-j} again we see that} $j(g_2z)$ satisfies the differential equation of $j(z)$.
\end{proof}

We consider a binary predicate $E^*_j(x,y)$ which will be interpreted in a differential field as 
$$ \exists y_1,y_2,y_3 \left( y_1^2y^2(y-1728)^2F(y,y_1,y_2,y_3)=0 \wedge y'=y_1x' \wedge y_1'=y_2x' \wedge y_2'=y_3x'\right).$$
Here we multiplied $F$ by $y_1^2y^2(y-1728)^2$ in order to make it a differential polynomial. Observe that any pair $(a,c)$, where $c$ is a constant, is in $E^*_j$. In order to simplify our arguments, we remove all points $(a,c)$ with $a \notin C, c \in C$, and define $E_j(x,y)$ by
$$E_j(x,y) \longleftrightarrow \left[ E^*_j(x,y) \wedge \neg (x'\neq 0 \wedge y'=0) \right] .$$
Actually, $E_j$ can be defined from $E^*_j$ (without using the derivation) as $C$ is definable by $E_j^*(0,y)$. The formula $E_j(0,y)$ defines the field of constants as well. 
One can also notice that for non-constant $x$ and $y$ the relation $E_j(x,y)$ is equivalent to $f(x,y)=0.$

\begin{definition}

The theory $T^0_j$ consists of the following first-order statements about a structure $K$ in the language $\mathfrak{L}_j:=\{ +, \cdot, E_j, 0, 1 \}$.

\begin{enumerate}
\item[A1] $K$ is an algebraically closed field. 

\item[A2] $C:=C_K=\{ c \in K: E_j(0,c)\}$ is an algebraically closed subfield. Further, $C^2 \subseteq E_j(K)$ and if $\left(z,j\right) \in E_j(K)$ and one of $z,j$ is constant then both of them are constants.

\item[A3] If $\left(z,j\right) \in E_j$ then for any  $g \in \SL_2(C)$, 
$\left(gz,j\right) \in E_j.$
Conversely, if for some $j$ we have $\left(z_1,j\right), \left(z_2,j\right) \in E_j$ then $z_2=gz_1$ for some $g \in \SL_2(C)$.

\item[A4] If $\left(z,j_1\right) \in E_j$ and $\Phi_N(j_1,j_2)=0$ for some $j_2$ and some modular polynomial $\Phi_N(X,Y)$ then $\left(z,j_2\right) \in E_j$.

\item[AS] If $\left(z_i,j_i\right) \in E_j,~ i=1,\ldots,n,$ with
$$\td_CC\left(\bar{z},\bar{j}\right) \leq n,$$ then $\Phi_N(j_i,j_k)=0$ for some $N$ and some $1\leq i < k \leq n$, or $j_i \in C$ for some $i$.
\end{enumerate}
\end{definition}

\begin{remark}\label{uniform-Ax}
A3 and A4 (the functional equations) imply that if $E_j(z_i,j_i),~ i=1,2$, and $j_1, j_2$ are modularly dependent then $z_1$ and $z_2$ have the same $\SL_2(C)$-orbit. However, the converse is not true: if $z_2=gz_1$ for some $g$ then this does not impose a relation on $j_1,j_2$ (they can be algebraically independent). Nevertheless, in that case we know by AS that $j_1$ and $j_2$ must be either algebraically independent or related by a modular relation (assuming $j_1$ and $j_2$ are non-constant).

A compactness argument shows that AS can be written as a first-order axiom scheme. Indeed, AS holds in all differential fields $K$. The compactness theorem can be applied to deduce that, given a parametric family of varieties $(W_{\bar{c}})_{\bar{c} \in C} \subseteq K^{2n}$, there is a natural number $N(W)$ such that if $\bar{c} \in C$ satisfies  $\dim W_{\bar{c}} \leq n$, and if $\left(\bar{z},\bar{j}\right) \in E_j(K)\cap W_{\bar{c}}(K)$ and $j_i \notin C$ for all $i$, then  $\Phi_N(j_i,j_k)=0$ for some $N \leq N(W)$ and some $1\leq i < k \leq n$. This can clearly be written as a first-order axiom scheme. Thus, AS should be understood as the uniform version of Ax-Schanuel.
\end{remark}

\begin{definition}
An $E_j$-\emph{field} is a model of $T^0_j$. If $K$ is an $E_j$-field, then a tuple $\left(\bar{z},\bar{j}\right) \in K^{2n}$ is called an $E_j$-\emph{point} if $\left(z_i, j_i\right) \in E_j(K)$ for each $i=1,\ldots,n$. By abuse of notation, we let $E_j(K)$ denote the set of all $E_j$-points in $K^{2n}$ for any natural number $n$ (which will be obvious from the context). The subfield $C_K$ is called the field of constants of $K$.
\end{definition}

The above lemmas show that reducts of differential fields to the language $\mathfrak{L}_j$ are $E_j$-fields. 

Let $C$ be an algebraically closed field with $\td(C/\mathbb{Q})=\aleph_0$ and let $\mathfrak{C}$ consist of all $E_j$-fields $K$ with $C_K=C$.  Note that $C$ is an $E_j$-field with $E_j(C)=C^2$ and it is the smallest structure in $\mathfrak{C}$. From now on, by an $E_j$-field we understand a member of $\mathfrak{C}$. Note that for some $X \subseteq A \in \mathfrak{C}$ we have $\langle X \rangle_A=C(X)^{\alg}$ (with the induced structure from $A$) and $\mathfrak{C}_{f.g.}$ consists of those $E_j$-fields that have finite transcendence degree over $C$. Therefore, a substructure of a finitely generated structure is finitely generated (Assumption \ref{fg}) and there is no infinite strictly descending chain of finitely generated structures (Assumption \ref{DCC}). 


\begin{definition}
For $A \subseteq B \in \mathfrak{C}_{f.g.}$ an $E_j$-\emph{basis} of $B$ over $A$ is an $E_j$-point $\bar{b}=\left(\bar{z},\bar{j}\right)$ from $B$ of maximal length satisfying the following conditions:
\begin{itemize}[leftmargin = 0.5cm]
\item $j_i$ and $j_k$ are modularly independent for all $i \neq k$,
\item $\left(z_i, j_i\right) \notin A^2$ for each $i$.
\end{itemize}
We let $\sigma(B/A)$ be the length of $\bar{j}$ in an $E_j$-basis of $B$ over $A$ (equivalently, $2\sigma(B/A)=| \bar{b}|$). When $A=C$ we write $\sigma(B)$ for $\sigma(B/C)$. Further, for $A \in \mathfrak{C}_{f.g.}$ define the predimension by 
$$\delta(A):=\td_{C}(A) - \sigma(A).$$
\end{definition}

Note that the Ax-Schanuel inequality for $j$ implies that $\sigma$ is finite for finitely generated structures. It is easy to see that for $A \subseteq B \in \mathfrak{C}_{f.g.}$ one has $\sigma(B/A)=\sigma(B)-\sigma(A)$. Moreover, for $A, B \subseteq D \in \mathfrak{C}_{f.g.}$ the inequality $$\sigma(AB) \geq \sigma(A)+\sigma(B)-\sigma(A \cap B)$$ holds. Hence $\delta$ is submodular (so it is a predimension) and the Pila-Tsimerman inequality states exactly that $\delta(A) \geq 0$ for all $A \in \mathfrak{C}_{f.g.}$ with equality holding if and only if $A=C$. The dimension associated with $\delta$ will be denoted by $d_j$ or simply $d$. We will add a superscript if we want to emphasise the model that we work in.

Observe also that for $A \subseteq B \in \mathfrak{C}_{f.g.}$
$$\delta(B/A)=\delta(B)-\delta(A)=\td(B/A)-\sigma(B/A).$$

\subsection{Amalgamation}\label{j-amalgamation}

\begin{definition}
A structure $A \in \mathfrak{C}$ is said to be \emph{full} if for every $j \in A$ there is $z \in A$ such that $A \models E_j\left(z,j\right)$. The subclass $\hat{\mathfrak{C}}$ consists of all full $E_j$-fields.
\end{definition}

\begin{lemma}\label{j-full-extension}
Every $A \in \mathfrak{C}$ has a unique (up to isomorphism over $A$) strong full extension $\hat{A} \in \hat{\mathfrak{C}}$ which is generated by $A$ as a full structure. In particular, if $A \in \mathfrak{C}_{f.g.}$ then $\hat{A} \in \hat{\mathfrak{C}}_{f.g.}$. Furthermore, if $f:A \hookrightarrow B$ is a strong embedding then $f$ extends to a strong embedding $\hat{f}:\hat{A} \hookrightarrow \hat{B}$.
\end{lemma}
\begin{proof}
Let $A \in \mathfrak{C}$. Choose an element $j \in A$ for which $A \models \neg \exists x E_j(x,j)$ (if there is such). Pick $z$ transcendental over $A$ (in a big algebraically closed field). Let $A_1:=A(z)^{\alg}$. Extend the relation $E_j$ to $A_1$ by adding the tuple $(z,j)$ to $E_j$ and closing the latter under the functional equations given by axioms A3 and A4. It is easy to see that $A \leq A_1$. Repeating this construction we will get a strong chain $A \leq A_1 \leq A_2 \leq \ldots$ the union of which, $A^1 := \bigcup_i A_i$, contains a solution of the formula $E_j(x,j)$ for each $j \in A$. Now we can iterate this construction and get another strong chain $A \leq A^1 \leq A^2 \leq \ldots$ such that for every $j \in A^i$ the formula $E_j(x,j)$ has a solution in $A^{i+1}$. The union $\hat{A}:=\bigcup_i A^i$ will be the desired strong and full extension of $A$. It is also clear that $\hat{A}$ is generated by $A$ as a full $E_j$-field.

Now we show that if $\hat{B} \in \hat{\mathfrak{C}}$ is a strong extension of $A$ then there is a strong embedding $\hat{A} \hookrightarrow \hat{B}$ over $A$. Let $j \in A$ be such that $A \models \neg \exists x E_j(x,j)$ and let $w\in \hat{B}$ satisfy $E_j(w,j)$. Since $w \notin A$, it must be transcendental over $A$. {We claim that $A_1$ (as constructed above) is isomorphic to $B_1:=A(w)^{\alg} \subseteq \hat{B}$ (with the induced structure). Indeed, $A \leq \hat{B}$ implies that $(w,j)$ is an $E_j$-basis of $B_1$ over $A$. Similarly, $(z,j)$ is an $E_j$-basis of $A_1$ over $A$. Hence, any isomorphism between the algebraically closed fields $A_1$ and $B_1$ that fixes $A$ pointwise and sends $z$ to $w$ is actually an isomorphism of $E_j$-fields $A_1$ and $B_1$.} Moreover, $B_1 \leq \hat{B}$ since $\delta(B_1/A)=0$. We can inductively construct similar partial isomorphisms from $\hat{A}$ into $\hat{B}$ the union of which will give a strong embedding $\hat{A} \hookrightarrow \hat{B}$. Furthermore, if $\hat{B}$ is generated by $A$ as a full $E_j$-field then we get an isomorphism $\hat{A} \cong \hat{B}$.
\end{proof}

\begin{proposition}\label{asym-amalgam}
The class $\hat{\mathfrak{C}}$ has the asymmetric amalgamation property.
\end{proposition}
\begin{proof}
Let $A, B_1, B_2 \in \hat{\mathfrak{C}}$ with an embedding $A \hookrightarrow B_1$ and a strong embedding $A \hookrightarrow B_2$. Let $B$ be the free amalgam of $B_1$ and $B_2$ over $A$ as algebraically closed fields. More precisely, $B$ is the algebraic closure of the extension of $A$ by the disjoint union of the transcendence bases of $B_1$ and $B_2$ over $A$. Identifying $A, B_1, B_2$ with their images in $B$ we have $B_1 \cap B_2 = A$.
Define $E_j$ on $B$ as the union $E_j(B_1) \cup E_j(B_2)$.

We show\footnote{Evidently, $B$ satisfies A1-A4 but we are still to prove that AS holds in $B$ too. So we do not know yet that $B$ is an $E_j$-field. However, $\delta$ is well defined on $B$ and it makes sense to say that $B_1 \leq B$.} that $B_1 \leq B$. By our definition of $E_j(B)$, a non-constant element $b \in B$ satisfies $B \models \exists x E_j(x,b)$ if and only if $b \in B_1 \cup B_2$. For a finitely generated $X \subseteq_{f.g.} B$ denote $X_1:=X \cap B_1, X_2:=X\cap B_2, X_0:=X \cap A$. From the above observation it follows that $\sigma(X)=\sigma(X \cap B_1)+ \sigma (X \cap B_2)-\sigma(X \cap A)$. Further, $X_1$ and $X_2$ are algebraically independent over $X_0$ and so 
$$\td(X/C)\geq \td(X_1/C)+\td(X_2/C)-\td(X_0/C).$$
Therefore
\begin{align*}
\delta(X)= & \td(X/C)-\sigma(X) \geq  \delta(X\cap B_1)+\delta(X \cap B_2)-\delta(X \cap A) \geq \delta(X \cap B_1),
\end{align*}
where the last inequality holds as $A \leq B_2$. Thus, $B_1 \leq B$.

This shows in particular that $\delta(X)\geq 0$. 
If $\delta(X)=0$ then $\delta(X \cap B_1)=0$ and so $X\cap B_1 \subseteq C$. But then $X \cap A \subseteq C$ which implies $\delta(X\cap A)=0$. Therefore $\delta(X\cap B_2)=0$ and $X \cap B_2 \subseteq C$. So $X \setminus C$ is disjoint from $B_1 \cup B_2$. But then $\delta(X)>0$ unless $X \subseteq C$.

So, $B$ satisfies the AS axiom scheme. Hence we can extend it strongly to a full $E_j$-field. The symmetric argument shows that if $A \leq B_1$ then $B_2 \leq B$.
\end{proof}

\begin{lemma}\label{countably-many-ext}
Let $A  \in \mathfrak{C}$ and let $B$ be a strong extension of $A$ finitely generated over $A$. Then $B$ is determined up to isomorphism by the locus $\Loc_A(\bar{b})$ for an $E_j$-basis $\bar{b}$ of $B$ over $A$ and the number $n=\td(B/A(\bar{b}))$. Hence for a given $A$ there are at most countably many strong finitely generated extensions of $A$, up to isomorphism.
\end{lemma}
\begin{proof}
{
Let $B_1$ and $B_2$ be two strong extensions of $A$, finitely generated over $A$. Let also $\bar{b}_i:= (\bar{z}_i,\bar{j}_i)$ be an $E_j$-basis of $B_i$ over $A$, and denote $A_i:=A(\bar{b}_i)$ ($i=1,2$). Assume $\Loc_A(\bar{b}_1) = \Loc_A(\bar{b}_2)$ and $\td(B/A_1) = \td(B/A_2)$. The map that fixes $A$ and sends $\bar{b}_1$ to $\bar{b}_2$ extends uniquely to a field isomorphism between $A_1$ and $A_2$, which respects the $E_j$-field structure. Any extension of this field isomorphism to $A_1^{\alg}$ and $A_2^{\alg}$ is actually an isomorphism of $E_j$-fields. Since $\bar{b}_i$ is an $E_j$-basis of $B_i$ over $A$, $E_j(B_i)=E_j(A_i^{\alg})$ for $i=1,2$. Therefore any extension of the above map to a field isomorphism of $B_1$ and $B_2$ (which exists as $\td(B/A_1) = \td(B/A_2)$) is an $E_j$-field isomorphism over $A$.}

For the second part of the lemma we just notice that there are countably many choices for $\Loc_A(\bar{b})$ and the number $n$.
\end{proof}

\begin{theorem}
The classes $\mathfrak{C}$ and $\hat{\mathfrak{C}}$ are strong amalgamation classes with the same strong Fra\"iss\'e limit $U$.
\end{theorem}
\begin{proof}
Proposition \ref{strong-amalgam} shows that $\mathfrak{C}$ has the strong amalgamation property. Then $\mathfrak{C}$ and $\hat{\mathfrak{C}}$ are strong amalgamation classes  and have the same strong Fra\"iss\'e limit by Proposition \ref{same-strong-limit}.
\end{proof}

Note that $\mathfrak{C}_{f.g.}$ does not have the asymmetric amalgamation property.

\subsection{Broad and free varieties}\label{broad-free}

\begin{definition}
Let $n$ be a positive integer, $k \leq n$ and $1\leq i_1 < \ldots < i_k \leq n$. Denote $\bar{i}=(i_1,\ldots,i_k)$ and define the projection map $\pr_{\bar{i}}:K^{n} \rightarrow K^{k}$ by
$$\pr_{\bar{i}}:(x_1,\ldots,x_n)\mapsto (x_{i_1},\ldots,x_{i_n}).$$
Further, define (by abuse of notation) $\pr_{\bar{i}}:K^{2n}\rightarrow K^{2k}$ by
$$\pr_{\bar{i}}:(\bar{x},\bar{y})\mapsto (\pr_{\bar{i}}\bar{x},\pr_{\bar{i}}\bar{y}).$$
\end{definition}

It will be clear from the context in which sense $\pr_{\bar{i}}$ should be understood (mostly in the second sense).

\begin{definition}
Let $K$ be an algebraically closed field. An irreducible algebraic variety $V \subseteq K^{2n}$ is \emph{{$E_j$-broad}} if and only if for any $0<k\leq n$ and any $1\leq i_1 < \ldots < i_k \leq n$ we have $\dim \pr_{\bar{i}} V \geq k$. We say $V$ is \emph{strongly $E_j$-broad} if the strict inequality $\dim \pr_{\bar{i}} V > k$ holds.\footnote{{In earlier versions of the paper and in \cite{Aslanyan-thesis} we called these varieties \emph{normal} following \cite{Zilb-pseudoexp}. The term \emph{broad} was coined by Sebastian Eterovi\'c.}}
\end{definition}

For a subfield $A \subseteq K$ we say a variety $V$ defined over $B \subseteq K$ is (strongly) $E_j$-broad over $A$ (regardless of whether $V$ is defined over $A$) if, for a generic (over $A \cup B$) point $\bar{v}$ of $V$, the locus $\Loc_A(\bar{v})$ is (strongly) $E_j$-broad. In other words, if a projection of $V$ is defined over $A$ then it is strongly $E_j$-broad over $A$.

\begin{remark}
It may seem strange that, in contrast to the exponential case, the functional equations of the $j$-function are not reflected in the definition of $E_j$-broadness. The reason is that those functional equations are of ``trivial'' type. Indeed, one would expect the following additional condition to be present: if $\bar{v}:=(\bar{z},\bar{j}) \in V(K)$ is a generic point of $V$ then we have not only $\td(\pr_{\bar{i}}\bar{v}) \geq k$ but also if we replace $z$'s by arbitrary elements in their $\SL_2(C)$-orbits and $j$'s by arbitrary elements in their Hecke orbits, then the transcendence degree of images of all those tuples under $\pr_{\bar{i}}$ must be at least $k$. However, it is obvious that the first condition already implies this because when we change the tuple in this manner, we do not change the transcendence degree (over $C$).
\end{remark}

\begin{remark}
$E_j$-broadness is a first-order definable property. This follows from the facts that irreducibility and algebraic dimension  are definable in algebraically closed fields. More generally, Morley rank is definable in strongly minimal theories.

Rotundity (in the exponential case) is first-order definable as well but it is not obvious since its definition involves infinitely many conditions which may be reduced to finitely many by weak CIT. See \cite[Theorem 4.9]{Kirby-semiab}.
\end{remark}

\begin{definition}
An algebraic variety $V \subseteq K^{2n}$ (with coordinates $\left(\bar{x},\bar{y}\right)$) is \emph{$E_j$-free} if it is not contained in any variety defined by an equation $\Phi_N(y_i,y_k)=0$ for a modular polynomial $\Phi_N$ and some indices $i, k$.
\end{definition}

This definition makes sense for an arbitrary field $K$. However, when $K$ is an $E_j$-field and $A \subseteq K$ is an $E_j$-subfield, we say $V \subseteq K^{2n}$ is \emph{$E_j$-free over} $A$ if it is $E_j$-free and it is not contained in a hyperplane defined by an equation of the form $y_i=a$ (for some $i$) where $a \in A$ with $A\models \exists z E_j(z,a)$.

We could require in the definition of $E_j$-freeness that $V$ is not contained in any variety defined by an equation of the form $y_i=b$ for some $b \in K$. This would be more standard definition and in fact it would be a definable property of the variety due to weak modular Zilber-Pink (but we will not need this result). Nevertheless, we find it more convenient to work with the notion of $E_j$-freeness (over $A$) defined above since it allows us to simplify some arguments slightly.

\begin{lemma}
If $A \leq B \in \mathfrak{C}_{f.g.}$ and $\bar{b} \in B^{2n}$ is an $E_j$-basis of $B$ over $A$ then the locus $\Loc_A(\bar{b})$ is $E_j$-broad and $E_j$-free over $A$, and strongly $E_j$-broad over $C$.
\end{lemma}
\begin{proof}
Follows obviously from definitions.
\end{proof}

\begin{lemma}\label{broad-strong-extension-j}
Let $A=C(\bar{a})^{\alg}$ be an $E_j$-field and $V$ be an $E_j$-broad irreducible variety defined over $A$. Then there is a strong extension $B$ of $A$ which contains an $E_j$-point of $V$ generic over $\bar{a}$. Furthermore, if $V$ is $E_j$-broad, $E_j$-free over $A$ and strongly $E_j$-broad over $C$ then we can choose $B$ so that $V(B)\cap E_j(B)$ contains a point generic in $V$ over $A$.
\end{lemma}
\begin{proof}
First, we prove the ``furthermore'' clause. Take a generic point of $V$ over $A$, say, $\bar{b}:=(\bar{z},\bar{j})$ and let $B:= \langle A \bar{b} \rangle = A(\bar{b})^{\alg}$. Extend $E_j$ by declaring $(z_i,j_i)$ an $E_j$-point for each $i$ and close it under functional equations (axioms A3 and A4). The given properties of $V$ make sure that $B$ is a model of $T^0_j$ and is a strong extension of $A$.

Now we prove the first part of the lemma. If for some $i_1< i_2< \ldots < i_k$ the projection $W:=\pr_{\bar{i}} V$ is defined over $C$ and has dimension $k$ then we pick constant elements $z_{i_s},j_{i_s}, s=1,\ldots,k,$ such that $(\bar{z},\bar{j})$ is generic in $W$ over $\bar{a}$. Doing this for all projections defined over $C$, we consider the variety $V_1$ obtained from $V$ by setting $x_{i_s}=z_{i_s},  y_{i_s}=j_{i_s}$ for all indices $i_s$ considered above. All of those pairs of constants will be in $E_j$.

Further, if $V_1$ is contained in a hyperplane $y_i=a$ for a non-constant $a \in A$ with $A \models E_j(z,a)$ for some $z \in A$, then we intersect it with the hyperplane $x_i=gz$ where we choose the entries of $g$ to be generic constants over $\bar{a}$. Doing this for all such $a$, we get a variety $V_2$, in a lower number of variables, which is still $E_j$-broad.

If $V_2$ is $E_j$-free then we proceed as above. Otherwise we argue as follows.
Suppose for some $i_1 \neq i_2$ the projection $\pr_{i_1,i_2} V_2$ satisfies the equation $\Phi_N(y_{i_1},y_{i_2})=0$ (we can assume $i_1$ and $i_2$ are different from all indices $i_s$ considered above). Let us assume for now that this is the only modular relation between the $y$-coordinates satisfied by $V_2$. Then we take algebraically independent elements $a, b, c \in C$ over $\bar{a}$ and over all elements from $A$ chosen above, and denote $d:=(1+bc)/a$. Let $V_3$ be the subvariety of $V_2$ defined by the equation $x_{i_2}=\frac{ax_{i_1}+b}{cx_{i_1}+d}$. It is easy to see that $\dim V_3 = \dim V_2 -1$ (here $V_3 \neq \emptyset$ as, by $E_j$-broadness, $\dim \pr_{i_1,i_2} V_2 \geq 2$). Now we take a generic point of $V_3$ over $\bar{a}abc$ and all constants taken above, and proceed as in the $E_j$-free case. Note that this generic point will be generic in $V$ over $A$.

{
When there are more modular relations between the $y$-coordinates of $V_2$, we apply the above procedure for all of those modular relations, that is, we introduce new generic $\SL_2(C)$-relations between the pairs of the appropriate $x$-coordinates (the corresponding $y$-coordinates of which satisfy a modular relation), and proceed as above.}
\end{proof}

\subsection{Existential closedness}\label{j-exist-clos}

Consider the following statements for an $E_j$-field $K$.


\begin{enumerate}

\item[EC] For each $E_j$-broad variety $V \subseteq K^{2n}$  the intersection $E_j(K) \cap V(K)$ is non-empty.

\item[SEC] For each $E_j$-broad variety $V \subseteq K^{2n}$ defined over a finite tuple  $\bar{a}\subseteq K$, the intersection $E_j(K) \cap V(K)$ contains a point generic in $V$ over $\bar{a}$.

\item[GSEC] For each irreducible variety $V \subseteq K^{2n}$ of dimension $n$ defined over a finitely generated strong $E_j$-subfield $A\leq K$, if $V$ is $E_j$-broad and $E_j$-free over $A$ and strongly $E_j$-broad over $C$, then the intersection $E_j(K) \cap V(K)$ contains a point generic in $V$ over $A$.


\item[NT] $K \supsetneq C$.

\item[ID] $K$ has infinite $d_j$-dimension.
\end{enumerate}

EC, SEC, GSEC, NT and ID stand for existential closedness, strong existential closedness, generic strong existential closedness,  non-triviality and infinite dimensionality respectively. Clearly, NT and EC are first-order axiomatisable. Notice that if an $E_j$-field $K$ satisfies AS+NT+EC then $\td(K/C)$ is infinite. In fact, all full $E_j$-fields with a non-constant point have the same property (we need to apply AS repeatedly).

\begin{lemma}\label{defined-over-C}
Let $V$ be an irreducible algebraic variety such that for every finitely generated (over $\mathbb{Q}$) field of definition $A \subseteq K$ there is a $C$-point generic in $V$ over $A$. Then $V$ is defined over $C$.
\end{lemma}
\begin{proof}
Let $A$ be a field of definition of $V$ and $\bar{a}$ be a transcendence basis of $A$ over $C$ (if $\bar{a}$ is empty then $V$ is defined over $C$). Then $V$ is defined over $\mathbb{Q}(\bar{a},\bar{c})^{\alg}$ for some finite tuple $\bar{c} \in C$. {Denote $A':=\mathbb{Q}(\bar{a},\bar{c})^{\alg}$. Let $\bar{d} \in C$ be a generic point of $V$ over $A$}. Then $\td(\bar{d}/A')=\dim V$. Since $\bar{a}$ is algebraically independent over $C$, we have $\td(\bar{d}/A')=\td(\bar{d}/\bar{c})$. Let $W:=\Loc(\bar{d}/C_0)$ where $C_0=\mathbb{Q}(\bar{c})^{\alg} \subseteq C$. Evidently, $W \supseteq V=\Loc(\bar{d}/A')$ and $\dim W = \td(\bar{d}/C_0)=\dim V$. Since both $V$ and $W$ are irreducible, $V=W$ and therefore $V$ is defined over $C_0$.
\end{proof}

\begin{proposition}\label{SEC-GSEC}
For $E_j$-fields  \emph{SEC} $\Rightarrow$ \emph{GSEC}.
\end{proposition}
\begin{proof}

Let $V$ and $A$ be as in the statement of GSEC. Choose $\bar{a} \subseteq A$ such that $\bar{a}$ contains an $E_j$-basis of $A$, $V$ is defined over $\bar{a}$ and $A=C(\bar{a})^{\alg}$.

Note that it suffices to prove that $V$ contains an $E_j$-point $\bar{v}=\left( \bar{z}, \bar{j} \right)$ none of the coordinates of which is constant and which is generic over $\bar{a}$. Indeed, we claim that $\bar{v}$ will be generic over $A$. If it is not the case then $\td(\bar{v}/A)<\dim V =n$. However, $j_i$ and $j_k$ are modularly independent for $i \neq k$ as $V$ is $E_j$-free and $\bar{v}$ is generic in $V$ over $\mathbb{Q}(\bar{a})$ and hence over $\mathbb{Q}$ (and modular polynomials are defined over $\mathbb{Q}$). Since $V$ is $E_j$-free over $A$ and $\bar{a}$ contains an $E_j$-basis of $A$, $\left(z_i,j_i\right) \notin A^2$ for each $i$. Then we would have $\delta(\bar{v}/A)<0$ which contradicts strongness of $A$ in $K$.

We claim that $V(K)$ contains an $E_j$-point generic over $\bar{a}$ which is not a $C$-point. If this is not the case then by SEC and Lemma \ref{defined-over-C} $V$ is defined over $C$. Since it is strongly $E_j$-broad over $C$, we have $\dim V >n+1$ which contradicts our assumption that $\dim V=n$. 

Now we prove that $V$ contains an $E_j$-point none of the coordinates of which is constant. We proceed to the proof by induction on $n$. The case $n=1$ is covered by the above argument (if $(z,j) \in E_j$ and one of $z,j$ is in $C$ then both of them must be in $C$). If $n>1$ take a point $\bar{v}=( \bar{z}, \bar{j} ) \in V(K)\cap E_j(K)$ generic over $\bar{a}$. If $\bar{v}$ has some constant coordinates then we can assume $(z_i,j_i) \subseteq C$ for $i=1,\ldots,k$ with $k<n$ (again, if one of $z_i,j_i$ is constant then both of them must be constants) and these are the only constant coordinates. If these constants have transcendence degree at least $k+1$ over $\bar{a}$ then the transcendence degree of all elements $z_i,j_i$ with $i>k$ over $C(\bar{a})$ will be strictly less than $n-k$ which contradicts $A \leq K$ as above.

Therefore $\td(\{ z_i,j_i:i\leq k \} / \mathbb{Q}(\bar{a}))=k$. By the induction hypothesis we can find an $E_j$-point $\bar{b}$ of $\pr_{\bar{i}}V$ (where $\bar{i}=(1,\ldots,k)$) none of the coordinates of which is constant and which is generic in $V$ over $\bar{a}$. Clearly, $\delta(\bar{b}/A)=0$ and so denoting $B:=A(\bar{b})^{\alg}$ we have $A \leq B \leq K$. Now let $V(\bar{b})$ be the variety obtained from $V$ by letting the corresponding $k$ coordinates of $V$ be equal to the corresponding coordinates of $\bar{b}$. Using the induction hypothesis we get an $E_j$-point $\bar{u}$ of $V(\bar{b})$ which is generic over $\bar{a},\bar{b}$ and whose coordinates are all non-constant. It is easy to see that $(\bar{b},\bar{u}) \in V(K)\cap E_j(K)$ is as required.
\end{proof}

\begin{proposition}
The strong Fra\"iss\'e limit $U$ satisfies \emph{SEC} and \emph{ID}, and hence \emph{GSEC}.
\end{proposition}
\begin{proof}
Let $V$ be an $E_j$-broad irreducible variety defined over a finite tuple $\bar{a}$. Let also $A:=\lceil \bar{a} \rceil_U$ (we can assume $A=C(\bar{a})^{\alg}$ by extending $\bar{a}$ if necessary). 
By Lemma \ref{broad-strong-extension-j} there is a strong extension $B$ of $A$ which contains an $E_j$-point $\bar{v}$ generic in $V$ over $\bar{a}$. Since $U$ is saturated for strong extensions, there is an embedding of $B$ into $U$ over $A$. The image of $\bar{v}$ under this embedding is the required generic $E_j$-point of $V$.  

For $n \in \mathbb{N}$, let $A_n$ be an algebraically closed field of transcendence degree $n$ over $C$. Defining $E_j(A_n)=C^2$ we make $A_n$ into a finitely generated $E_j$-field with $d_j$-dimension $n$. By universality of $U$, $A_n$ can be strongly embedded into $U$ which shows $U$ has infinite $d_j$-dimension because strong extensions preserve dimension.
\end{proof}

One can directly prove in the same manner that $U$ satisfies GSEC (without using Proposition \ref{SEC-GSEC}).

\begin{lemma}\label{generic-hyperplanes}
Let $K$ be an infinite $d_j$-dimensional $E_j$-field and $A\subseteq K$ be a finitely generated $E_j$-subfield. Assume $V \subseteq K^{2n}$ is an $E_j$-broad irreducible variety defined over $A$ with $\dim V > n$. Then we can find a strong extension $A \leq A' \leq K$, generated over $A$ by finitely many $d_j^K$-independent (over $A$) elements, and an $E_j$-broad subvariety $V'$ of $V$, defined over $A'$, with $\dim V'=n$.

\end{lemma}

This can be proven exactly as in the exponential case by intersecting $V$ with generic hyperplanes (see \cite{Kirby-semiab}, Proposition 2.33 and Theorem 2.35). We give full details for completeness.

For $\bar{p}:=(p_1, \ldots, p_N) \in K^N \setminus \{ 0 \}$ let the hyperplane $\Pi_{\bar{p}}$ be defined by the equation $ \sum_{i=1}^N p_ix_i=1.$ It is obvious that $\bar{a} \in \Pi_{\bar{b}}(K)$ iff $\bar{b} \in \Pi_{\bar{a}}(K)$.

We will need the following result from \cite{Kirby-semiab} (Lemma 2.31) which has been adapted from \cite{Zilber-bicolored}. 

\begin{lemma}\label{hyperplane-lemma}
Let $A\subseteq K$ be an algebraically closed field,  $\bar{v} \in K^N$ and  $\bar{p} \in \Pi_{\bar{v}}$ be generic over $A(\bar{v})$. Then for any tuple $\bar{w} \in A(\bar{v})^{\alg}$ either $\bar{v} \in A(\bar{w})^{\alg}$ or $\td(\bar{w}/A\bar{p})=\td(\bar{w}/A)$ (i.e. $\bar{w} \forksym_A \bar{p}$ in $\ACF_0$).
\end{lemma}
\begin{proof}
{Assume $\bar{w} \nforkind_A \bar{p}$. Then $\td(\bar{p}/A\bar{w})<\td(\bar{p}/A)$. We may also assume $\bar{v}\nsubseteq A$ as otherwise the result is trivial. In this case $\td(\bar{p}/A)=N$. 

Let $P:=\Loc(\bar{p}/A(\bar{w})^{\alg})$. Then $$\dim P =\td(\bar{p}/A\bar{w})\geq \td(\bar{p}/A\bar{v})=N-1,$$
the inequality following from the fact that $\bar{w} \in A(\bar{v})^{\alg}$. On the other hand, $\td(\bar{p}/A\bar{w})<\td(\bar{p}/A) = N$. Thus, $\dim P = N-1$.

Observe that  $P \supseteq \Pi_{\bar{v}}$, for $\Pi_{\bar{v}}=\Loc(\bar{p}/A\bar{v})$. Since both $P$ and $\Pi_{\bar{v}}$ are irreducible and have the same dimension, they must be equal.} Hence $\Pi_{\bar{v}}$ is defined over $A(\bar{w})^{\alg}$ and so the formula $\forall \bar{y} \in \Pi_{\bar{v}} (\bar{x} \in \Pi_{\bar{y}})$ defines $\bar{v}$ over $A(\bar{w})^{\alg}$. Thus, $\bar{v} \in A(\bar{w})^{\alg}$.
\end{proof}

\begin{proof}[Proof of Lemma \ref{generic-hyperplanes}]
Let $\dim V> n$. It will be enough to find $A'$ and $V'$ with $\dim V'=\dim V-1$. Pick a generic point $\bar{v} \in V(K)$. Denote $N=2n$ and choose $p_1,\ldots,p_{N-1} \in K$ to be $d_j^K$-independent over $A$. Pick $p_N \in K$ such that $\sum_{i=1}^N p_i v_i =1$.

Let $A':=A(\bar{p})^{\alg}$ and $V':=V \cap \Pi_{\bar{p}}=\Loc_{A'}(\bar{v})$. Obviously, $V'$ is irreducible and $\dim V' = \dim V-1$. We claim that $V'$ is $E_j$-broad. Let $\bar{w}:=\pr_{\bar{i}} \bar{v}$ for some projection map $\pr_{\bar{i}}$ with $| \bar{i} |=k\leq n$. Then obviously $\bar{w} \in A(\bar{v})^{\alg}$. Therefore by Lemma \ref{hyperplane-lemma} either $\bar{v} \in A(\bar{w})^{\alg}$ or $\td(\bar{w}/A')=\td(\bar{w}/A)$. In the former case $$\dim \pr_{\bar{i}} V'= \td(\bar{w}/A')=\td(\bar{v}/A')=\dim V' = \dim V -1 \geq n \geq k.$$
In the latter case $$\dim \pr_{\bar{i}} V'= \td(\bar{w}/A')=\td(\bar{w}/A)\geq k,$$
where the last inequality follows from $E_j$-broadness of $V$.
\end{proof}

\begin{proposition}
The strong Fra\"iss\'e limit $U$ is the unique countable $E_j$-field satisfying \emph{GSEC} and \emph{ID} and having $\td(C/\mathbb{Q})=\aleph_0$.
\end{proposition}
\begin{proof}
Let $K$ be such an $E_j$-field. We will show it is saturated with respect to strong embeddings. Let $A \leq B$ be finitely generated $E_j$-fields and let $\bar{b}$ be a basis of $B$ over $A$. If $\td(B/A(\bar{b}))>0$ then let $\bar{b}'$ be a transcendence basis of $B$ over $A(\bar{b})$. We can find a strong extension $B \leq B'=B(\bar{a}')^{\alg}$ such that $B' \models  E_j(a'_i,b'_i)$ for each $i$. Replacing $B$ by $B'$ we may assume that $\td(B/A(\bar{b}))=0$ and hence $B=A(\bar{b})^{\alg}$.

Let $V:=\Loc_A(\bar{b})$ be the Zariski closure of $\bar{b}$ over $A$. It is irreducible, $E_j$-broad and $E_j$-free over $A$ and strongly $E_j$-broad over $C$. By Lemma \ref{generic-hyperplanes} we can find a strong extension $A'$ of $A$, generated by independent elements over $A$, and an $E_j$-broad irreducible subvariety $V'$ of $V$ over $A'$ such that $\dim V'= \sigma(B/A)$. Obviously, $V'$ is also $E_j$-free over $A'$ and strongly $E_j$-broad over $C$ (because $V$ is).

By GSEC there is a point $\bar{v} \in V' \cap E_j$ in $K$, generic in $V'$ over $A'$. Then $\bar{v}$ is also generic in $V$ over $A$. Let $B'':=A'(\bar{v})^{\alg}$ with the induced structure from $K$. Then $\delta(A')=\delta(B'')$ and so $B'' \leq K$. Now $B':=A(\bar{v})^{\alg}$ with the induced structure is isomorphic to $B$ over $A$. Moreover, $B''$ is generated by $d_j$-independent elements over $B'$ and so $B' \leq B''$ and $B' \leq K$. Therefore, $K$ is saturated for strong extensions.
\end{proof}

\subsection{The complete theory}

\begin{definition}
Let $T_j$ be the theory axiomatised by $T^0_j+ \mbox{EC}+\mbox{NT}$.
\end{definition}

Note that $T_j$ is an $\forall \exists$-theory. 

\begin{proposition}\label{saturated-SEC}
All $\aleph_0$-saturated models of $T_j$ satisfy \emph{SEC}, and hence \emph{GSEC}. 
\end{proposition}
\begin{proof}
It suffices to show that in an arbitrary model $K$ of $T_j$ every Zariski-open subset of an irreducible $E_j$-broad variety contains an $E_j$-point. Let $(\bar{x},\bar{y})$ be the coordinates of $K^{2n}$ and let $V \subseteq K^{2n}$ be  an $E_j$-broad irreducible variety. It is enough to show that for every proper subvariety $W$ of $V$, defined by a single equation, $V\setminus W$ contains an $E_j$-point. Suppose $W$ (as a subvariety of $V$) is defined by an equation $f(z_1,\ldots,z_k)=0$ where each $z_i$ is one of the coordinates $\{ x_i, y_i: i=1,\ldots,n \}$. The assumption that $W \subsetneq V$ means that $f$ does not vanish on $V$. 

We use Rabinovich's trick to replace $V\setminus W$ by an $E_j$-broad irreducible variety in  a higher number of variables. Consider the variety $V' \subseteq K^{2(n+1)}$ (with coordinates $(\bar{x},x_{n+1},\bar{y},y_{n+1})$) defined by the equations of $V$ and one additional equation $x_{n+1}f(\bar{z})=1$. It is clear that $V'$ is $E_j$-broad and irreducible. By EC, $V'$ contains an $E_j$-point. Its projection onto the coordinates  $(\bar{x},\bar{y})$ will be an $E_j$-point in $V\setminus W$.
\end{proof}

\begin{proposition}\label{ID}
All $\aleph_0$-saturated models of $T_j$ satisfy \emph{ID}. In particular, a countable saturated model of $T_j$ (if it exists) is isomorphic to $U$.
\end{proposition}

We will need the following algebraic lemma in the proof.

\begin{lemma}\label{lemma-id}
Let $(K;+,\cdot, C)$ be an $\aleph_0$-saturated pair of algebraically closed fields and $k$ and $n$ be a positive integer. Then there is a strongly $E_j$-broad variety $P\subseteq K^{2n}$ defined over $C$ of dimension $n+1$ such that for any $\bar{c}\in C^k$ the variety $P$ is not contained in a proper subvariety of $K^{2n}$ defined over $\mathbb{Q}(\bar{c})^{\alg}$.
\end{lemma}
\begin{proof}
Consider a variety $P$ defined by $n-1$ equations of the form $y_i=f_i(\bar{c}_i,\bar{x})$ where $f_i(\bar{c}_i,\bar{X})$ is a polynomial of $\bar{X}$ over $\bar{c}_i$. If we choose the degrees of $f_i$'s (hence the lengths of $\bar{c}_i$'s) sufficiently large, and letting all  $\bar{c}_i$'s be algebraically independent over $\mathbb{Q}$, then $P$ will be as desired.
\end{proof}

\begin{proof}[Proof of Proposition \ref{ID}]
Let $K \models T_j$ be $\aleph_0$-saturated. A priori, we do not have a type whose realisations would be $d_j$-independent, but we can write $d_j$-independence by an $\mathfrak{L}_{\omega_1,\omega}$-sentence. The idea is to use weak modular Zilber-Pink to reduce this $\mathfrak{L}_{\omega_1,\omega}$-sentence to a type and show that it is finitely satisfiable in $K$.

ID means that for each $n$ there is a $2n$-tuple $\bar{x}$ of algebraically independent (over $C$) elements with $\bar{x} \in E_j(K)$ (which is equivalent to $\delta(\bar{x})=n$) such that for all tuples $\bar{y}$ one has $\delta(\bar{y},\bar{x})\geq n$. Here we can assume as well that $\bar{y}$ is a $2l$-tuple for some $l$ and is an $E_j$-point. The fact that $\bar{x}$ is algebraically independent over $C$ is given by a type consisting of formulae $\varphi_i(\bar{x})= \forall \bar{c} (\bar{x} \notin V_i(\bar{c})),~ i<\omega,$ stating that $\bar{x}$ is not in any hypersurface (defined over $C$) from a parametric family of hypersurfaces $(V_i(\bar{c}))_{\bar{c}\in C}$ (to be more precise, we could say that $(V_i(\bar{c}))_{\bar{c} \in C}$ is the parametric family of hypersurfaces over $C$ of degree $i$). 

The statement $\forall y_1, \ldots, y_{2l} \delta(\bar{y},\bar{x})\geq n$ can be written as an $\mathfrak{L}_{\omega_1,\omega}$-sentence as follows. Given an algebraic variety $W \subseteq K^{2l+2n+m}$ defined over $\mathbb{Q}$, for any $\bar{c} \in C^m$ with $\dim W(\bar{c})<2n+l$ and for any $\bar{y} \in W(\bar{x},\bar{c})\cap E_j$, the $j$-coordinates of $\bar{y}$ (i.e. $y_{l+1},\ldots,y_{2l}$) must satisfy a modular relation $\Phi_N(y_{l+i},y_{l+k})=0$ for some $N$ and some $1\leq i<k \leq l$, or a modular relation with $\bar{x}$, i.e. $\Phi_N(x_{n+i},y_{l+k})=0$ for some $1\leq i \leq n,~ 1\leq k \leq l$, or we must have $y_{l+i} \in C$ for some $1\leq i \leq l$.

Now suppose, for contradiction, that ID does not hold in $K$. It means that for some $n$, for all $2n$-tuples $\bar{x}$ satisfying $\bar{x} \in E_j$ and $\bigwedge_i \varphi_i(\bar{x})$, there are a variety $W \subseteq K^{2l+2n+m}$ (for some $l,m$) defined over $\mathbb{Q}$, a constant point $\bar{c} \in C^m$ with $\dim W(\bar{c})<2n+l$, and a tuple $\bar{y} \in W(\bar{x},\bar{c})\cap E_j$, such that $\Phi_N(y_{l+i},y_{l+k})\neq 0$ for all $N$ and all $1\leq i<k \leq l$, and $\Phi_N(x_{n+i},y_{l+k})\neq 0$ for all $1\leq i \leq n,~ 1\leq k \leq l$, and $y_{l+i} \notin C$ for all $1\leq i \leq l$.


For a parametric family of varieties $W(\bar{c})_{\bar{c} \in C^m}$ in $K^{2l+2n}$ let $N(W)$ be the maximal number $N$ such that $\Phi_N$ occurs in the defining equations of the finitely many special varieties given by the weak modular Zilber-Pink for this parametric family. 
Then the following holds\footnote{It holds for any $N$ instead of $N(W)$. Our choice of $N(W)$ was made so that it will lead to a contradiction.} in $K$:
\begin{gather*}
\forall \bar{x}  \Bigg[  \bar{x}\in E_j \wedge \bigwedge_{i< \omega} \varphi_i(\bar{x}) \longrightarrow \bigvee_{\substack{l,m \in \mathbb{N} \\ W \subseteq K^{2l+2n+m}}} \exists \bar{c} \in C^m ~\exists \bar{y} \in W(\bar{x},\bar{c}) \cap E_j  
\\
\Bigg( \dim W(\bar{c})<2n+l~ 
 \wedge \bigwedge_{\substack{p\leq N(W)\\ 1\leq i<k \leq l}} \Phi_p(y_{l+i},y_{l+k}) \neq 0 \\ \wedge \bigwedge_{\substack{p\leq N(W)\\ 1\leq i \leq n \\ 1 \leq k \leq l}} \Phi_p(x_{n+i},y_{l+k}) \neq 0 ~
 \wedge \bigwedge_{1\leq i \leq l} y_{l+i} \notin C \Bigg) \Bigg].
\end{gather*}
Here the disjunction (in the first line) is over all positive integers $l, m$ and all algebraic varieties $W \subseteq K^{2l+2n+m}$ defined over $\mathbb{Q}$ (there are countably many such triples $(l,m,W)$).

By $\aleph_0$-saturation of $K$ and compactness we deduce that there are a finite collection of varieties $W_s \subseteq K^{2l_s+2n+m_s},~ s=1,\ldots,t$, and a finite number $r$ such that

\begin{gather*}
\forall \bar{x}  \Bigg[  \bar{x}\in E_j \wedge \bigwedge_{i\leq r} \varphi_i(\bar{x}) \longrightarrow \bigvee_{s\leq t} \exists \bar{c} \in C^{m_s} ~\exists \bar{y} \in W_s(\bar{x},\bar{c}) \cap E_j \\ 
\Bigg( \dim W_s(\bar{c})<2n+l_s ~
  \wedge  \bigwedge_{\substack{p\leq N(W_s)\\ 1\leq i<k \leq l_s}} \Phi_p(y_{l_s+i},y_{l_s+k}) \neq 0 \\
  \wedge \bigwedge_{\substack{p\leq N(W_s)\\ 1\leq i \leq n \\ 1 \leq k \leq l_s}} \Phi_p(x_{n+i},y_{l_s+k}) \neq 0  ~\wedge \bigwedge_{1\leq i \leq l_s} y_{l_s+i} \notin C \Bigg) \Bigg].
\end{gather*}

The formulas $\varphi_i(\bar{u}) = \forall \bar{c} (\bar{u} \notin V_i(\bar{c}))$ state that $\bar{u}$ is not in a given parametric family of hypersurfaces $V_i(\bar{c})$. It is easy to see that we can find a strongly $E_j$-broad and $E_j$-free variety $P$ in $K^{2n}$ defined over $C$, of dimension $n+1$, which is not contained in any of the varieties $V_i(\bar{c})$ for any $\bar{c}$ and any $i\leq r$.  Moreover, we can choose $P$ as in Lemma \ref{lemma-id} with  $k = \max \{ m_s: s\leq t \}$.

Now by the GSEC property we can find a non-constant point $\bar{a} = (a_1,\ldots,a_{2n}) \in E_j(K) \cap P(K)$ which is generic in $P$ over $C$. Indeed, we need to intersect $P$ with a generic hyperplane as in Lemma \ref{generic-hyperplanes}, with algebraically independent coefficients (instead of $d_j$-independent), and get an $E_j$-broad and $E_j$-free variety over (the strong closure of) the field generated by those coefficients. Then we apply GSEC. By Lemma \ref{lemma-id} the elements $a_1,\ldots,a_{2n}$ are algebraically independent over any subfield $C_0 \subseteq C$ with $\td(C_0/\mathbb{Q})\leq k$.


Then $\td(\bar{a}/C)=n+1$ and $a_{n+1},\ldots,a_{2n}$ are pairwise modularly independent (since $P$ is $E_j$-free), hence $\delta(\bar{a})=1$. Moreover, $\varphi_i(\bar{a})$ holds for $i \leq r$. Therefore by the above statement, for some $W:=W_s \subseteq K^{2l+2n+m}$ (where $m=m_s,l=l_s$) there are $\bar{c} \in C^m, \bar{b} \in W(\bar{a},\bar{c})(K) \cap E_j(K)$ such that $\dim W(\bar{c})<2n+l$ and $a_{n+1},\ldots,a_{2n},b_{l+1},\ldots,b_{2l}$ are non-constant and do not satisfy any modular equation $\Phi_p=0$ for $p \leq N(W)$. By our choice of $k$ we also know that $\td(\bar{a}/\bar{c})=2n$.

Suppose, for a moment, that $a_{n+1},\ldots,a_{2n},b_{l+1},\ldots,b_{2l}$ are pairwise modularly independent. Then evidently $\delta(\bar{a},\bar{b})\leq 0$ which contradicts AS.

However those elements may satisfy some modular relations $\Phi_p=0$ with $p>N(W)$. Let $S\subseteq K^{2l+2n}$ be the special variety defined by all those modular relations (more precisely, $S$ is a component of the variety defined by those relations which contains the point $(\bar{b},\bar{a})$).\footnote{Let us stress  again that $S$ is defined only by the modular relations satisfied by the tuple $(\bar{b}', \bar{a}') =   (b_{l+1},\ldots,b_{2l}, d_{n+1},\ldots,d_{2n})$. In particular, there are modular relations only between $n+l$ coordinates of $S$.} Let also $R \subseteq S \cap W(\bar{c})$ be a component of the intersection containing that point. We claim that $S$ intersects $W(\bar{c})$ typically, i.e. $R$ is a typical component of the intersection (in $K^{2l+2n}$). Indeed, by our choice of $N(W)$, the intersection cannot be strongly atypical. On the other hand, no coordinate is constant on $R$ since $b_i,d_k \notin C$, so $R$ is not an atypical component. Here we actually need to show that $R$ does not satisfy any equation of the form $r_i= d$ where $r_i$ is the $i$-th coordinate of $R$ and $d$ is a fixed element of $K$ (and not necessarily of $C$). However, since $R$ is defined over $C$, such an element $d$ would necessarily be from $C$.

It means that if $a_{n+1},\ldots,a_{2n},b_{l+1},\ldots,b_{2l}$ satisfy $h$ independent modular relations (i.e. $h=\codim S = 2n+2l-\dim S$), then
$$\dim R = \dim S + \dim W(\bar{c}) - (2l+2n) =  \dim W(\bar{c})-h.$$

Since $a_1, \ldots, a_{2n}$ are algebraically independent over $\mathbb{Q}(\bar{c})$ and $R$ is defined over $\mathbb{Q}(\bar{c})^{\alg}$, we have
$$\dim R(\bar{a}) = \dim R - 2n = \dim W(\bar{c})-2n - h < l-h.$$

Then we have 

\begin{gather*}
     \td(\bar{b}/C(\bar{a})) \leq \td(\bar{b}/\bar{a},\bar{c}) \leq \dim R(\bar{a}) \leq l-h-1.
\end{gather*}

Thus,

\begin{gather*}
    \td(\bar{b}, \bar{a}/C) = \td(\bar{b}/C(\bar{a})) + \td(\bar{a}/C) \leq (l-h-1) + (n+1) = n+l-h.
\end{gather*}

On the other hand $\sigma (\bar{a},\bar{b}) = n+l-h$. So $\delta(\bar{a},\bar{b})=0$ which contradicts AS.
\end{proof}

\begin{proposition}
The theory $T_j$ is complete and the Fra\"iss\'e limit $U$ is $\aleph_0$-saturated.
\end{proposition}
\begin{proof}
{Let $T_j^1$ be an arbitrary completion of $T_j$ and let $M$ be a (possibly uncountable) $\aleph_0$-saturated model of $T_j^1$. Let also $C:=C_M$ be the field of constants (which may be uncountable as well). 
}


\begin{claim}
For all finitely generated (i.e. of finite transcendence degree over $C$) strong $E_j$-subfields $A, B \leq M$ with an isomorphism $f:A \cong B$, and for any $a' \in M$, there are $A \leq A' \leq M$ and $B\leq B' \leq M$ with $a' \in A'$ such that $f$ extends to an isomorphism $A'\cong B'$.
\end{claim}
\begin{proof}[Proof of the claim]

We can assume $a' \notin A$ and hence it is transcendental over $A$. We consider two cases.



Case 1: $d_j^M(a'/A)=0$.

Let $A':=\lceil Aa' \rceil_M$ and let $\bar{v}$ be an $E_j$-basis of $A'$ over $A$. Since $\delta(A'/A)=0$, $A'$ must be $\mathfrak{C}$-generated by $\bar{v}$ over $A$, i.e. $A'=\langle A\bar{v} \rangle$. Now if $V:=\Loc_A(\bar{v})$, then $V$ is $E_j$-broad and $E_j$-free over $A$ and strongly $E_j$-broad over $C$, and $\dim V =n$ (since $\delta(\bar{v}/A)=0$). Let $W$ be the image of $V$ under the isomorphism $f : A \rightarrow B$ (i.e. we just replace the coefficients of equations of $V$ by their images under $f$). Then $W$ is $E_j$-broad and $E_j$-free over $B$ and strongly $E_j$-broad over $C$, and so by the GSEC property the intersection $W(M) \cap E_j(M)$ contains a point $\bar{w}$ generic in $W$ over $B$. Setting $B':=B(\bar{w})^{\alg}$ (with the induced structure from $M$), we see that $\delta(B'/B)=0$ and so $B \leq B' \leq M$. Clearly $f$ extends to an isomorphism from $A'$ to $B'$.

Case 2: $d_j^M(a'/A)=1$.

In this case we pick an element $b'\in M$ which is $d_j^M$-independent from $B$ (which exists by ID) and set $A'=\langle Aa' \rangle = A(a')^{\alg}$ and $B'=\langle Bb'\rangle = B(b')^{\alg}$. Obviously $A' \leq M,~ B' \leq M$ and $A' \cong B'$.
\end{proof}

{Thus, given two tuples $\bar{a},\bar{b} \in M$ (of the same length) with an isomorphism $f:\lceil \bar{a} \rceil_M \cong \lceil \bar{b} \rceil_M$ sending $\bar{a}$ to $\bar{b}$,  we can start with $f$ and construct a back-and-forth system of partial isomorphisms from $M$ to itself showing that $\tp^M(\bar{a}) = \tp^M(\bar{b})$. Combining this with Lemma \ref{countably-many-ext} we see that if $A:=\lceil \bar{a} \rceil_M$ and $\bar{a}'$ is an $E_j$-basis of $A$ then the type of $\bar{a}$ in $M$ is determined uniquely by $\Loc(\bar{a}'/C)$, $\Loc(\bar{a}/C(\bar{a}'))$ and the number $\td(A/C(\bar{a}'))$. Indeed, if for $\bar{a},\bar{b} \in M$ these data coincide then there is an isomorphism $f:\lceil \bar{a} \rceil_M \cong \lceil \bar{b} \rceil_M$ sending $\bar{a}$ to $\bar{b}$. Moreover, if $\Loc(\bar{a}'/C)$ and $\Loc(\bar{a},\bar{a}'/C)$ are defined over a finite set of constants $\bar{c}$, then the proof of Lemma \ref{countably-many-ext} shows actually that $\tp^M(\bar{a})$ is determined by the algebraic varieties $\Loc(\bar{c}/\mathbb{Q})$, $\Loc(\bar{a}',\bar{c}/\mathbb{Q})$ and $\Loc(\bar{a},\bar{a}',\bar{c}/\mathbb{Q})$ (in fact, the first two varieties are also uniquely determined by the third one) and the number $\td(A/C(\bar{a}'))$. There are countably many choices for those varieties and the transcendence degree, hence $T^1_j$ is small, i.e. there are countably many pure types (types over $\emptyset$).\footnote{We can in fact show by a similar argument that $T_j$ is $\aleph_0$-stable.} This implies that $T^1_j$ has a countable saturated model which must be isomorphic to $U$ by Proposition \ref{ID}. Thus, $U$ is saturated and $T_j^1 = \Th(U)$. Since $T_j^1$ was an arbitrary completion of $T_j$, the latter has a unique completion and so it is complete.}
\end{proof}

We summarise the results of this section in the following theorems.

\begin{theorem}\label{thm: j axiomatisation}
The theory $T_j$ is consistent and complete. It is the first-order theory of the strong Fra\"iss\'e limit $U$, which is saturated.
\end{theorem}

\begin{theorem}\label{thm: j main}
The following are equivalent.
\begin{itemize}[leftmargin = 0.5cm]
\item The Ax-Schanuel inequality for $j$ is adequate.
\item The Ax-Schanuel inequality for $j$ is strongly adequate.
\item $\mathfrak{L}_{j}$-reducts of differentially closed fields are models of $T_j$.
\item $\mathfrak{L}_{j}$-reducts of differentially closed fields satisfy \emph{EC}.
\item $\mathfrak{L}_{j}$-reducts of $\aleph_0$-saturated differentially closed fields satisfy \emph{SEC}.
\end{itemize}
\end{theorem}

Thus, adequacy of the Ax-Schanuel inequality for $j$ gives a complete axiomatisation of the first-order theory of the differential equation of $j$ and show that it is nearly model complete. It also gives a criterion for a system of differential equations in terms of the equation of $j$ to have a solution. 

It was recently shown in \cite{Aslanyan-Eterovic-Kirby-Diff-EC-j} that differentially closed fields satisfy EC for $j$. 

\begin{theorem}[{\cite[Theorem 1.1]{Aslanyan-Eterovic-Kirby-Diff-EC-j}}]\label{fact:ec-j}
$\mathfrak{L}_{j}$-reducts of differentially closed fields satisfy \emph{EC}.
\end{theorem}

Combining Theorems \ref{thm: j main} and \ref{fact:ec-j} we get the following result.

\begin{corollary}
The Ax-Schanuel inequality for $j$ is strongly adequate, and $T_j$ is the complete theory of $\mathfrak{L}_{j}$-reducts of differentially closed fields.
\end{corollary}

\section{The $j$-function with derivatives}\label{j-general-case}

In this section we study the predimension given by the Ax-Schanuel inequality \emph{with derivatives}.
We consider a predicate $E'_j(x,y,y_1,y_2)$ which will be interpreted in a differential field as 
$$ \exists y_3 \left( y_1^2y^2(y-1728)^2F(y,y_1,y_2,y_3)=0 \wedge y'=y_1x' \wedge y_1'=y_2x' \wedge y_2'=y_3x'\right).$$
Note that all quadruples of constants $(z,j,j^{(1)},j^{(2)})$ satisfy $E'_j$ unless $j^{(1)}=0,~ j^{(2)} \neq 0$. For convenience we extend $E'_j$ so that it contains all quadruples of constants. Also, if $z$ is constant then $j,j^{(1)},j^{(2)}$ must be constants as well. Moreover, if $\bar{a}=\left(z_i,j_i,j^{(1)}_i,j^{(2)}_i\right) \in E'_j(K^{\times})$ and one of the coordinates of $\bar{a}$ is constant then all of them are. One can also notice that for non-constant $x$ and $y$ the relation $E'_j$ is equivalent to $$f(x,y)=0\wedge y_1=\partial_x y \wedge y_2=\partial^2_x y.$$

\subsection{$E_j'$-fields}

\begin{lemma}
In a differential field if $x_2=gx_1$ with $g=\begin{pmatrix}
a & b \\
c & d
\end{pmatrix}
\in \SL_2(C)$ then for any non-constant $y$ we have
\begin{gather*}
\partial_{x_2}y=\partial_{x_1}y \cdot (cx_1+d)^2,\\ 
\partial^2_{x_2}y=\partial^2_{x_1}y \cdot (cx_1+d)^2 -2c \cdot \partial_{x_1}y\cdot (cx_1+d)^3.
\end{gather*}
\end{lemma}
\begin{proof}
Easy calculations.
\end{proof}

\begin{definition}

The theory $(T^0_j)'$ consists of the following first-order statements about a structure $K$ in the language $\mathfrak{L}_j:=\{ +, \cdot, E'_j, 0, 1 \}$.

\begin{enumerate}
\item[A1'] $K$ is an algebraically closed field with an algebraically closed subfield $C:=C_K$, which is defined by $E'_j(0,y,0,0)$. Further, $C^4 \subseteq E'_j(K)$ and if $\bar{a}=\left(z,j,j^{(1)},j^{(2)}\right) \in E'_j(K^{\times})$ and one of the coordinates of $\bar{a}$ is in $C$ then $\bar{a} \subseteq C$.

\item[A2'] For any $z, j \in K \setminus C$ there is at most one pair $\left(j^{(1)},j^{(2)}\right)$ in $K$ with $E'_j\left(z,j,j^{(1)},j^{(2)}\right)$.

\item[A3'] If $\left(z,j,j^{(1)},j^{(2)}\right) \in E'_j$ then for any  $g=\begin{pmatrix} a & b \\ c & d \end{pmatrix} \in \SL_2(C)$ 
$$\left(gz,j,j^{(1)}\cdot (cz+d)^2,j^{(2)}\cdot (cz+d)^2-2c\cdot j^{(1)}\cdot (cz+d)^3\right) \in E'_j.$$
Conversely, if for some $j$ we have $\left(z_1,j,j^{(1)},j^{(2)}\right), \left(z_2,j, w^{(1)},w^{(2)}\right) \in E'_j$ then $z_2=gz_1$ for some $g \in \SL_2(C)$.

\item[A4'] If $\left(z,j_1,j_1^{(1)},j_1^{(2)}\right) \in E'_j$ and $\Phi(j_1,j_2)=0$ for some modular polynomial $\Phi(X,Y)$ then $\left(z,j_2,j_2^{(1)},j_2^{(2)}\right) \in E'_j$ where $j_2^{(1)},~ j_2^{(2)}$ are determined from the following system of equations:
\begin{align*}
& \frac{\partial \Phi}{\partial X}(j_1,j_2)\cdot j_1^{(1)} + \frac{\partial \Phi}{\partial Y}(j_1,j_2)\cdot j_2^{(1)} =0,\\
& \frac{\partial^2 \Phi}{\partial X^2}(j_1,j_2)\cdot \left(j_1^{(1)}\right)^2+ \frac{\partial^2 \Phi}{\partial Y^2}(j_1,j_2)\cdot \left(j_2^{(1)}\right)^2 + 2\frac{\partial^2 \Phi}{\partial X \partial Y}(j_1,j_2)\cdot j_1^{(1)}\cdot j_2^{(1)}+ \\ 
& \frac{\partial \Phi}{\partial X}(j_1,j_2)\cdot j_1^{(2)} + \frac{\partial \Phi}{\partial Y}(j_1,j_2)\cdot j_2^{(2)}=0.
\end{align*}

\item[AS'] If $\left(z_i,j_i,j_i^{(1)},j_i^{(2)}\right) \in E'_j,~ i=1,\ldots,n,$ with
$$\td_CC\left(\bar{z},\bar{j},\bar{j}^{(1)},\bar{j}^{(2)}\right) \leq 3n$$ then $\Phi_N(j_i,j_k)=0$ for some $N$ and $1\leq i < k \leq n$ or $j_i \in C$ for some $i$.
\end{enumerate}
\end{definition}

A4' is obtained by differentiating the equality $\Phi(j_1,j_2)=0$. A compactness argument shows that AS' can be written as a first-order axiom scheme exactly as before.

\begin{definition}
An $E'_j$-\emph{field} is a model of $(T^0_j)'$. If $K$ is an $E'_j$-field, then a tuple $\left(\bar{z},\bar{j},\bar{j}^{(1)},\bar{j}^{(2)}\right) \in K^{4n}$ is called an $E'_j$-\emph{point} if $\left(z_i, j_i\right) \in E'_j(K)$ for each $i=1,\ldots,n$. By abuse of notation, we let $E'_j(K)$ denote the set of all $E'_j$-points in $K^{4n}$ for any natural number $n$.
\end{definition}

Let $C$ be an algebraically closed field with $\td(C/\mathbb{Q})=\aleph_0$ and let $\mathfrak{C}$ consist of all $E'_j$-fields $K$ with $C_K=C$.  Note that $C$ is an $E'_j$-field with $E'_j(C)=C^4$ and it is the smallest structure in $\mathfrak{C}$. From now on, by an $E'_j$-field we understand a member of $\mathfrak{C}$. 


\begin{definition}
For $A \subseteq B \in \mathfrak{C}_{f.g.}$ an $E'_j$-\emph{basis} of $B$ over $A$ is an $E'_j$-point $\bar{b}=\left(\bar{z},\bar{j},\bar{j}^{(1)},\bar{j}^{(2)}\right)$ from $B$ of maximal length satisfying the following conditions:
\begin{itemize}
\item $j_i$ and $j_k$ are modularly independent for all $i \neq k$,
\item $\left(z_i, j_i, j^{(1)},j^{(2)}_i\right) \notin A^4$ for each $i$.
\end{itemize}
We let $\sigma(B/A)$ be the length of $\bar{j}$ in an $E'_j$-basis of $B$ over $A$ (equivalently, $4\sigma(B/A)=| \bar{b}|$). When $A=C$ we write $\sigma(B)$ for $\sigma(B/C)$. Further, for $A \in \mathfrak{C}_{f.g.}$ define the predimension by 
$$\delta(A):=\td_{C}(A) - 3\cdot\sigma(A).$$
\end{definition}

As before, $\delta$ is submodular (so it is a predimension) and the Pila-Tsimerman inequality states exactly that $\delta(A) \geq 0$ for all $A \in \mathfrak{C}_{f.g.}$ with equality holding if and only if $A=C$. The dimension associated with $\delta$ will be denoted by $d_j'$.


\begin{definition}
A structure $A \in \mathfrak{C}$ is said to be \emph{full} if for every $j \in A$ there are $z,j^{(1)},j^{(2)} \in A$ such that $A \models E'_j\left(z,j,j^{(1)},j^{(2)}\right)$. The subclass $\hat{\mathfrak{C}}$ consists of all full $E'_j$-fields.
\end{definition}

The obvious analogues of all results from Sections \ref{j-univtheory-predim} and \ref{j-amalgamation} hold in this setting as well (with obvious adaptations of the proofs). So we get a strong Fra\"iss\'e limit $U$.

\subsection{Existential Closedness}

\begin{definition}
Let $n$ be a positive integer, $k \leq n$ and $1\leq i_1 < \ldots < i_k \leq n$. Denote $i=(i_1,\ldots,i_k)$ and define the projection map $\pr_{\bar{i}}:K^{n} \rightarrow K^{k}$ by
$$\pr_{\bar{i}}:(x_1,\ldots,x_n)\mapsto (x_{i_1},\ldots,x_{i_n}).$$
Further, define $\pr_{\bar{i}}:K^{4n}\rightarrow K^{4k}$ by
$$\pr_{\bar{i}}:(\bar{x},\bar{y},\bar{z},\bar{w})\mapsto (\pr_{\bar{i}}\bar{x},\pr_{\bar{i}}\bar{y},\pr_{\bar{i}}\bar{z},\pr_{\bar{i}}\bar{w}).$$
\end{definition}

Below $\pr_{\bar{i}}$ should always be understood in the second sense.

\begin{definition}
Let $K$ be an algebraically closed field. An irreducible algebraic variety $V \subseteq K^{4n}$ is \emph{$E'_j$-broad} if and only if for any $1\leq i_1 < \ldots < i_k \leq n$ we have $\dim \pr_{\bar{i}} V \geq 3k$. We say $V$ is \emph{strongly $E'_j$-broad} if the strict inequality $\dim \pr_{\bar{i}} V > 3k$ holds.
\end{definition}



\begin{definition}
An algebraic variety $V \subseteq K^{4n}$ (with coordinates $\left(\bar{x},\bar{y}, \bar{y}^{(1)}, \bar{y}^{(2)}\right)$) is \emph{$E'_j$-free} if it is not contained in any variety defined by an equation $\Phi_N(y_i,y_k)=0$ for some modular polynomial $\Phi_N$ and some indices $i, k$.

When $K$ is an $E'_j$-field and $A \subseteq K$ is an $E'_j$-subfield, we say $V \subseteq K^{4n}$ is \emph{$E'_j$-free over} $A$ if it is $E'_j$-free and it is not contained in a hyperplane defined by an equation $y_i=a$ (for some $i$) where $a \in A$ with $A\models \exists z, u, v E'_j(z,a,u,v)$.
\end{definition}

Consider the following statements for an $E'_j$-field $K$.


\begin{enumerate}

\item[EC'] For each $E'_j$-broad variety $V \subseteq K^{4n}$  the intersection $E'_j(K) \cap V(K)$ is non-empty.

\item[SEC'] For each $E'_j$-broad variety $V \subseteq K^{4n}$ defined over a finite tuple  $\bar{a}\subseteq K$, the intersection $E'_j(K) \cap V(K)$ contains a point generic in $V$ over $\bar{a}$.

\item[GSEC'] For each irreducible variety $V \subseteq K^{4n}$ of dimension $3n$ defined over a finitely generated strong $E'_j$-subfield $A\leq K$, if $V$ is $E'_j$-broad and $E'_j$-free over $A$ and strongly $E'_j$-broad over $C$, then the intersection $E'_j(K) \cap V(K)$ contains a point generic in $V$ over $A$.


\item[NT'] $K \supsetneq C$.

\item[ID'] $K$ has infinite $d_j'$-dimension.
\end{enumerate}

Again, the analogues of all facts established in Sections \ref{broad-free} and \ref{j-exist-clos} are true with obvious adaptations of the proofs. Therefore $U$ is a model of the theory $T'_j$ axiomatised by A1-A4,AS,NT,EC. 

Proposition \ref{saturated-SEC} holds as well. However, in order to prove that $\aleph_0$-saturated models of $T_j'$ satisfy ID' (Proposition \ref{ID}) we need to use a Zilber-Pink type statement for $j$ and its derivatives. This significantly complicates the situation, and we discuss the details in the next section.

{

\subsection{Functional Modular Zilber-Pink with Derivatives}

In this subsection we state the Functional Modular Zilber-Pink with Derivatives (FMZPD) theorem, which we need to prove infinite dimensionality. See \cite{Aslanyan-weakMZPD} for details on FMZPD.

In this section special varieties defined in Definition \ref{def-special} will be called $j$-special.

\begin{definition}
Let $C$ be an algebraically closed field. A $C$-\emph{geodesic variety} $U\subseteq C^n$ (with coordinates $\bar{x}$) is an irreducible component of a variety defined by equations of the form $x_i=g_{i,k}x_k$ for some $g_{i,k}\in \SL_2(C)$. If $S\subseteq C^n$ (with coordinates $\bar{y}$) is a $j$-special variety, then $U$ is said to be a $C$-geodesic variety associated with $S$ if for any $1\leq i, k\leq n$ we have $\Phi_N(y_i,y_k)=0$ on $S$ for some $N$ if and only if $x_i=g_{i,k}x_k$ on $U$ for some $g_{i,k}\in \SL_2(C)$.
\end{definition}

\begin{definition}
Let $C$ be an algebraically closed field. Define $D$ as the zero derivation on $C$ and extend $(C;+,\cdot, D)$ to a differentially closed field $(K;+,\cdot, D)$. 
\begin{itemize}[leftmargin = 0.5cm]
    \item Let $T \subseteq C^n$ be a $j$-special variety and $U\subseteq C^n$ be a $C$-geodesic variety associated with $T$. Denote by $\langle \langle U, T \rangle \rangle$ the Zariski closure over $C$ of the projection of the set $$(E_j')^{\times}(K) \cap (U(K)\times T(K) \times K^2)$$ onto the last $3n$ coordinates where $(E_j')^{\times}$ is the set of all $E_j'$-points with no constant coordinates.
    
    \item A D-\emph{special} variety is a variety $S:=\langle \langle U, T \rangle \rangle$ for some $T$ and $U$ as above. In this case $S$ is said to be a D-special variety associated with $T$ and $U$. We will also say that $T$ (or $U$) is a $j$-special (respectively, geodesic) variety associated with $S$. A D-special variety associated with $T$ is one associated with $T$ and $U$ for some $C$-geodesic variety $U$ associated with $T$. 
    
    
\end{itemize}
\end{definition}

\begin{definition}
Let $T\subseteq C^n$ be a $j$-special variety. A geodesic variety $U$ associated with $T$ is called \emph{upper triangular} if all matrices $g$ occurring in the definition of $U$ are upper triangular. If $U$ is upper triangular then a D-special variety associated with $T$ and $U$ is also called \emph{upper triangular}. 
\end{definition}

Note that if $S$ is an upper triangular D-special variety associated with $T$ then $\dim S = 3 \dim T$.

\begin{definition}
Let $V\subseteq C^{3n}$ be an algebraic variety (or, more generally, an arbitrary set). \emph{A D-special closure} of $V$ is a D-special variety $S\subseteq C^{3n}$ which contains $V$ and is minimal among the D-special varieties containing $V$. 
\end{definition}

By Noetherianity of the Zariski topology every variety has at least one D-special closure which, in general, is not unique. 

\begin{definition}\label{Def-D-broad}
Let $V \subseteq C^{3n}$ be an irreducible variety with a D-special closure $S$ and let $T$ be the $j$-special closure of $\pr_{\bar{y}} V$. Then $V$ is said to be D-\emph{broad} if for all $1\leq i_1 < \ldots < i_k \leq n$ 
\begin{equation*}
    \dim \pr_{\bar{i}} V \geq \dim \pr_{\bar{i}} S - \dim \pr_{\bar{i}} T.
\end{equation*}
If for all $\bar{i}$ the above inequality is strict then $V$ is \emph{strongly D-broad}.
\end{definition}

Strong D-broadness of a variety does not depend on the choice of its D-special closure.

\begin{definition}
Let $V\subseteq C^{3n}$ be a variety. A \textit{D-atypical} subvariety of $V$ in $C^{3n}$ is an atypical component $W$ of an intersection $V\cap T$ where $T\subseteq C^{3n}$ is D-special. If, in addition, $W$ is strongly D-broad then we say that it is \textit{strongly D-atypical}.
\end{definition}

\begin{theorem}[FMZPD]\label{thm-fmzpd}
Let $C$ be an algebraically closed field. Given a parametric family of algebraic varieties $(V_{\bar{c}})_{\bar{c} \subseteq C}$ of $C^{3n}$, there is a finite collection $\Sigma$ of proper $j$-special subvarieties of $C^n$ such that for every $\bar{c} \subseteq C$, every strongly D-atypical subvariety of $V_{\bar{c}}$ is contained in a D-special variety associated with some $T\in \Sigma$.
\end{theorem}

We will need a slight generalisation of this theorem where we work in $C^{4n}$ rather than $C^{3n}$. A D-special subvariety of $C^{4n}$ is a variety of the form $C^n\times S$ where $S\subseteq C^{3n}$ is D-special. A component of an intersection of $V$ with a D-special subvariety is strongly D-atypical if it is an atypical component of the intersection and it is strongly D-broad (as in Definition \ref{Def-D-broad}).

\begin{corollary}\label{horizontal-FMZPD}
Let $C$ be an algebraically closed field. Given a parametric family of algebraic varieties $(V_{\bar{c}})_{\bar{c} \subseteq C}$ of $C^{4n}$, there is a finite collection $\Sigma$ of proper $j$-special subvarieties of $C^n$ such that for every $\bar{c} \subseteq C$, every strongly D-atypical subvariety of $V_{\bar{c}}$ is contained in a D-special variety associated with some $S\in \Sigma$.
\end{corollary}
\begin{proof}
This can be deduced from Theorem \ref{thm-fmzpd} as in the proof of \cite[Theorem 11.4]{Bays-Kirby-exp}. Indeed, the statement follows from Theorem \ref{thm-fmzpd} and the fibre dimension theorem if all fibres of the projection of $V_{\bar{c}}$ on the first $n$ coordinates (i.e. the $z$-coordinates) have the same dimension. If there are fibres with different dimensions then we can consider the definable subsets of a fixed dimension (which has finitely many possible values) and apply the theorem to each of those separately.
\end{proof}

\subsection{Infinite dimensionality}

\begin{proposition}\label{ID'}
All $\aleph_0$-saturated models of $T_j'$ satisfy \emph{ID'}. In particular, a countable saturated model of $T'_j$ (if it exists) is isomorphic to $U$.
\end{proposition}
\begin{proof}
This is a generalisation of the proof of Proposition \ref{ID}.
 
Let $K \models T'_j$ be $\aleph_0$-saturated. ID means that for each $n$ there is a $4n$-tuple $\bar{x}$ of algebraically independent (over $C$) elements with $\bar{x} \in E'_j(K)$ (which is equivalent to $\delta(\bar{x})=n$) such that for all tuples $\bar{y}$ one has $\delta(\bar{y},\bar{x})\geq n$. Here we can assume as well that $\bar{y}$ is a $4l$-tuple for some $l$ and is an $E'_j$-point. The fact that $\bar{x}$ is algebraically independent over $C$ is given by a type consisting of formulae $\varphi_i(\bar{x})= \forall \bar{c} (\bar{x} \notin V_i(\bar{c})),~ i<\omega,$ stating that $\bar{x}$ is not in any hypersurface (defined over $C$) from a parametric family of hypersurfaces $(V_i(\bar{c}))_{\bar{c}\in C}$ (here $(V_i(\bar{c}))_{\bar{c} \in C}$ is the parametric family of hypersurfaces in $K^{4n}$ over $C$ of degree $i$). 

The statement $\forall y_1, \ldots, y_{4l} \delta(\bar{y},\bar{x})\geq n$ can be written as an $\mathfrak{L}_{\omega_1,\omega}$-sentence as follows. Given an algebraic variety $W \subseteq K^{4l+4n+m}$ defined over $\mathbb{Q}$, for any $\bar{c} \in C^m$ with $\dim W(\bar{c})<4n+3l$ and for any $\bar{y} \in W(\bar{x},\bar{c})\cap E'_j$, the $j$-coordinates of $\bar{y}$ (i.e. $y_{l+1},\ldots,y_{2l}$) must satisfy a modular relation $\Phi_N(y_{l+i},y_{l+k})=0$ for some $N$ and some $1\leq i<k \leq l$, or a modular relation with $\bar{x}$, i.e. $\Phi_N(x_{n+i},y_{l+k})=0$ for some $1\leq i \leq n,~ 1\leq k \leq l$, or we must have $y_{l+i} \in C$ for some $1\leq i \leq l$.

Now suppose, for contradiction, that ID' does not hold in $K$. It means that for some $n$, for all $4n$-tuples $\bar{x}$ satisfying $\bar{x} \in E'_j$ and $\bigwedge_i \varphi_i(\bar{x})$, there are a variety $W \subseteq K^{4l+4n+m}$ (for some $l,m$) defined over $\mathbb{Q}$, a constant point $\bar{c} \in C^m$ with $\dim W(\bar{c})<4n+3l$, and a tuple $\bar{y} \in W(\bar{x},\bar{c})\cap E'_j$, such that $\Phi_N(y_{l+i},y_{l+k})\neq 0$ for all $N$ and all $1\leq i<k \leq l$, and $\Phi_N(x_{n+i},y_{l+k})\neq 0$ for all $1\leq i \leq n,~ 1\leq k \leq l$, and $y_{l+i} \notin C$ for all $1\leq i \leq l$.

For a parametric family of varieties $W(\bar{c})_{\bar{c} \in C^m}$ in $K^{4l+4n}$ let $N(W)$ be the maximal number $N$ such that $\Phi_N$ occurs in the defining equations of the D-special varieties given by Corollary \ref{horizontal-FMZPD}. Then the following holds in $K$:
\begin{gather*}
\forall \bar{x}  \Bigg[  \bar{x}\in E_j \wedge \bigwedge_{i< \omega} \varphi_i(\bar{x}) \longrightarrow \bigvee_{\substack{l,m \in \mathbb{N} \\ W \subseteq K^{4l+4n+m}}} \exists \bar{c} \in C^m ~\exists \bar{y} \in W(\bar{x},\bar{c}) \cap E_j  
\\
\Bigg( \dim W(\bar{c})<4n+3l~ 
 \wedge \bigwedge_{\substack{p\leq N(W)\\ 1\leq i<k \leq l}} \Phi_p(y_{l+i},y_{l+k}) \neq 0 \\ \wedge \bigwedge_{\substack{p\leq N(W)\\ 1\leq i \leq n \\ 1 \leq k \leq l}} \Phi_p(x_{n+i},y_{l+k}) \neq 0 ~
 \wedge \bigwedge_{1\leq i \leq l} y_{l+i} \notin C \Bigg) \Bigg].
\end{gather*}
Here the disjunction (in the first line) is over all positive integers $l, m$ and all algebraic varieties $W \subseteq K^{2l+2n+m}$ defined over $\mathbb{Q}$ (there are countably many such triples $(l,m,W)$).

By $\aleph_0$-saturation of $K$ and compactness we deduce that there are a finite collection of varieties $W_s \subseteq K^{4l_s+4n+m_s},~ s=1,\ldots,t$, and a finite number $r$ such that

\begin{gather*}
\forall \bar{x}  \Bigg[  \bar{x}\in E_j \wedge \bigwedge_{i\leq r} \varphi_i(\bar{x}) \longrightarrow \bigvee_{s\leq t} \exists \bar{c} \in C^{m_s} ~\exists \bar{y} \in W_s(\bar{x},\bar{c}) \cap E_j \\ 
\Bigg( \dim W_s(\bar{c})<4n+3l_s ~
  \wedge  \bigwedge_{\substack{p\leq N(W_s)\\ 1\leq i<k \leq l_s}} \Phi_p(y_{l_s+i},y_{l_s+k}) \neq 0 \\
  \wedge \bigwedge_{\substack{p\leq N(W_s)\\ 1\leq i \leq n \\ 1 \leq k \leq l_s}} \Phi_p(x_{n+i},y_{l_s+k}) \neq 0  ~\wedge \bigwedge_{1\leq i \leq l_s} y_{l_s+i} \notin C \Bigg) \Bigg].
\end{gather*}

As in the proof of Proposition \ref{ID}, we can find a strongly $E'_j$-broad and $E'_j$-free variety $P$ in $K^{4n}$ defined over $C$, of dimension $3n+1$, which is not contained in any of the varieties $V_i(\bar{c})$ for any $\bar{c}$ and any $i\leq r$ (recall that $(V_i(\bar{c}))_{\bar{c} \in C}$ is the parametric family of hypersurfaces over $C$ of degree $i$).

Now by the GSEC property we can find a non-constant point $\bar{a} = (a_1,\ldots,a_{4n}) \in E_j'(K)\cap P(K)$ which is generic in $P$ over $C$. Moreover, as in Lemma \ref{lemma-id}, we could choose $P$ ``generic'' enough so that the elements $a_1,\ldots,a_{2n}$ are algebraically independent over any subfield $C_0 \subseteq C$ with $\td(C_0/\mathbb{Q})\leq k$ where $k := \max \{ m_s: s\leq t \}$.

Then $\td(\bar{a}/C)=3n+1$ and $a_{n+1},\ldots,a_{2n}$ are pairwise modularly independent, hence $\delta(\bar{a})=1$. Moreover, $\varphi_i(\bar{a})$ holds for $i \leq r$. Therefore by the above statement, for some $W:=W_s \subseteq K^{4l+4n+m}$ (where $m=m_s,l=l_s$) there are $\bar{c} \in C^m, \bar{b} \in W(\bar{a},\bar{c})(K) \cap E_j(K)$ such that $\dim W(\bar{c})<4n+3l$ and $a_{n+1},\ldots,a_{2n},b_{l+1},\ldots,b_{2l}$ are non-constant and do not satisfy any modular equation $\Phi_p=0$ for $p \leq N(W)$. By changing some of the coordinates of $\bar{b}$ if necessary, we may assume that all $\GL_2(C)$-relations among $a_1,\ldots,a_n,b_1,\ldots,b_l$ are given by the identity matrix. This transformation will change the variety $W(\bar{c})$ sending it to its image under the action of some elements of $\GL_2(C)$. Since that way we still get a parametric family of varieties and $\dim W(\bar{c})$ does not change, this will not cause any issues. By our choice of $k$ we also know that $\td(\bar{a}/\bar{c})=4n$.


The tuple $(\bar{b},\bar{a})$ may satisfy some modular relations $\Phi_p=0$ with $p>N(W)$. Let $T\subseteq K^{3l+3n}$ be the D-special closure of the projection of $(\bar{b},\bar{a})$ onto the last $3l$ coordinates of $\bar{b}$ and last $3n$ coordinates of $\bar{a}$. Denote $S:=K^{n+l}\times T$. Assume $R \subseteq S \cap W(\bar{c})$ is a component of the intersection containing the point $(\bar{b},\bar{a})$. We claim that $S$ intersects $W(\bar{c})$ typically, i.e. $R$ is a typical component of the intersection (in $K^{4l+4n}$). Indeed, by our choice of $N(W)$, the intersection cannot be strongly D-atypical. On the other hand, since $R$ contains an $E_j'$-point none of the coordinates of which are in $C$, by Ax-Schanuel it must be strongly D-broad.
So $R$ is not an atypical component.

It means that if $a_{n+1},\ldots,a_{2n},b_{l+1},\ldots,b_{2l}$ satisfy $h$ independent modular relations, then $\dim T = 3(l+n-h),~ \dim S = 4(n+l)-3h$ (for we assumed $T$ is upper triangular), and
$$\dim R = \dim S + \dim W(\bar{c}) - (4l+4n) =  \dim W(\bar{c})-3h.$$

Since $a_1, \ldots, a_{4n}$ are algebraically independent over $\mathbb{Q}(\bar{c})$ and $R$ is defined over $\mathbb{Q}(\bar{c})^{\alg}$, we have
$$\dim R(\bar{a}) = \dim R - 4n = \dim W(\bar{c})-4n - 3h < 3(l-h).$$

Then we have 

\begin{gather*}
     \td(\bar{b}/C(\bar{a})) \leq \td(\bar{b}/\bar{a},\bar{c}) \leq \dim R(\bar{a}) \leq 3(l-h)-1.
\end{gather*}

Thus,

\begin{gather*}
    \td(\bar{b}, \bar{a}/C) = \td(\bar{b}/C(\bar{a})) + \td(\bar{a}/C) \leq 3(l-h)-1 + (3n+1) = 3(n+l-h).
\end{gather*}

On the other hand $\sigma (\bar{a},\bar{b}) = n+l-h$. So $\delta(\bar{a},\bar{b})=0$ which contradicts AS.
\end{proof}}

\subsection{The complete theory}

Now, as in Section \ref{j-chapter}, we obtain the following results.

\begin{theorem}\label{thm: j' axiomatisation}
The theory $T'_j$ is consistent and complete. It is the first-order theory of the strong Fra\"iss\'e limit $U$, which is saturated.
\end{theorem}

\begin{theorem}\label{thm: j' main}
The following are equivalent.
\begin{itemize}[leftmargin = 0.5cm]
\item The Ax-Schanuel inequality for $j$ and its derivatives is adequate.
\item The Ax-Schanuel inequality for $j$ and its derivatives is strongly adequate.
\item $\mathfrak{L}_{j}'$-reducts of differentially closed fields are models of $T_j'$.
\item $\mathfrak{L}'_{j}$-reducts of differentially closed fields satisfy \emph{EC'}.
\item $\mathfrak{L}'_{j}$-reducts of $\aleph_0$-saturated differentially closed fields satisfy \emph{SEC'}.
\end{itemize}
\end{theorem}

Like EC for $E_j$, EC for $E_j'$ was also established in \cite{Aslanyan-Eterovic-Kirby-Diff-EC-j}.

\begin{theorem}[{\cite[Theorem 1.2]{Aslanyan-Eterovic-Kirby-Diff-EC-j}}]\label{fact:ec-j'}
$\mathfrak{L}'_{j}$-reducts of differentially closed fields satisfy \emph{EC'}.
\end{theorem}

\begin{corollary}
The Ax-Schanuel inequality for $j$ and its derivatives is strongly adequate, and $T_j'$ is the complete theory of $\mathfrak{L}'_{j}$-reducts of differentially closed fields.
\end{corollary}

\addtocontents{toc}{\protect\setcounter{tocdepth}{1}}
\subsection*{Acknowledgements} I am grateful to Boris Zilber, Jonathan Pila, Jonathan Kirby, Ehud Hrushovski, Felix Weitkamper and Sebastian Eterovi\'c  for numerous useful discussions and comments. I would also like to thank the referee for many valuable comments on the paper.

\addtocontents{toc}{\protect\setcounter{tocdepth}{2}}

\bibliographystyle {alpha}
\bibliography {ref}

\end{document}